\documentclass[11pt]{amsart}
\usepackage[T1]{fontenc}
\usepackage[utf8]{inputenc}
\usepackage{amsmath,amsthm,amssymb}
\usepackage{newtxtext,newtxmath}
\usepackage[a4paper,margin=1.1in]{geometry}
\usepackage{booktabs}
\usepackage{microtype}
\usepackage[hidelinks]{hyperref}
\theoremstyle{plain}
\newtheorem{theorem}{Theorem}[section]
\newtheorem{corollary}[theorem]{Corollary}
\newtheorem{lemma}[theorem]{Lemma}
\newtheorem{proposition}[theorem]{Proposition}
\theoremstyle{definition}
\newtheorem{definition}[theorem]{Definition}
\newtheorem{example}[theorem]{Example}
\theoremstyle{remark}
\newtheorem{remark}[theorem]{Remark}
\numberwithin{equation}{section}
\newcommand{\FPP}{\operatorname{FPP}}

\newcommand{\Prob}{\mathbb{P}}
\newcommand{\GL}{\operatorname{GL}}
\newcommand{\Deck}{\operatorname{Deck}}
\newcommand{\N}{\mathbb{N}}
\newcommand{\Aut}{\operatorname{Aut}}
\newcommand{\End}{\operatorname{End}}
\newcommand{\Gal}{\operatorname{Gal}}
\newcommand{\Fix}{\operatorname{Fix}}

\newcommand{\PP}{\mathbf{P}}
\newcommand{\OO}{\mathcal{O}}
\newcommand{\Z}{\mathbb{Z}}
\newcommand{\Q}{\mathbb{Q}}
\newcommand{\C}{\mathbb{C}}
\newcommand{\Fp}{\mathbb{F}}
\newcommand{\abs}[1]{\lvert #1\rvert}

\title[Fixed-point proportion of Latt\`es maps]{The fixed-point proportion of Latt\`es maps}
\author{Xin Li}
\address{Shenzhen University of Advanced Technology, Shenzhen, China}
\email{suat25000420@stu.suat-sz.edu.cn}
\date{November 5, 2026}
\hypersetup{pdftitle={The fixed-point proportion of Latt\`es maps},pdfauthor={Xin Li}}
\begin{document}
\begin{abstract}
How often does an element of the geometric iterated Galois group of a Latt\`es map fix a vertex of the tree? Let $A/K$ be an abelian variety over a separably closed field of characteristic $0$, $G\le\Aut(A,0)$ finite, and $\alpha$ an isogeny of degree at least $2$ normalising $G$. The Latt\`es map on $X=A/G$ has a preimage tree governed by $\alpha$-power torsion; we determine the fixed-point proportion in every dimension. For $\zeta\in G$ the weight is the reciprocal of the number of fixed points of $\zeta$ in $\alpha$-power torsion (zero if infinite); equivalently, it is a product of $p$-adic absolute values of determinants on the $\alpha$-active subspaces; for scalar $\alpha=m$ the proportion is $|G|^{-1}\sum_{1\ne\zeta\in G}\prod_{p\mid m}\abs{\det(1-\zeta\mid H_1(A,\Z))}_p$. We deduce positivity and vanishing criteria, prove the groups $G_\infty$ are martingale groups, determine the limit law, and reduce groups with translations to the fixed-origin case. In dimension one we determine ten values for scalar multipliers and nontrivial $G\le\Aut(E,0)$, and for the eight quadratic Latt\`es classes the proportions are $1/8$, $1/4$ and $5/16$. Bridy--Jones--Kelsey--Lodge exclude the Latt\`es case from their main theorem; their Theorem~2.5 periodic-point bound is attained, via reductions modulo split primes $p\equiv1\pmod4$, for the map $(x^2+a)/(x^2-a)$, $a^2=-1$, and we settle the exceptional case of the Fari\~na-Asategui--Radi conjecture.
\end{abstract}
\maketitle

\section{Introduction}\label{sec:intro}

A natural question about an arboreal Galois group is how often its elements fix a vertex of the tree. Its fixed-point proportion --- the limiting proportion of elements that fix at least one vertex at every level --- measures exactly this, and for Latt\`es maps the preimage tree is governed by torsion, which makes the invariant computable. We determine the fixed-point proportion in every dimension. The Latt\`es case is one of the two cases excluded from the main theorem (Theorem~1.3) of Bridy--Jones--Kelsey--Lodge: for a class of quadratic Latt\`es maps they prove a $1/8$ lower bound for the periodic-point proportion \cite[Theorem~2.5]{BJKL}, and for the Latt\`es map of Corollary~\ref{cor:BJKL} this bound is attained.

Let $K$ be a separably closed field of characteristic $0$, $A/K$ an abelian variety, $G\le\Aut(A,0)$ a finite group, and $\alpha\in\End(A)$ an isogeny of degree at least $2$ with $\alpha G\alpha^{-1}=G$. Then $\alpha$ descends to a finite map $f\colon X\to X$ on $X=A/G$. For a point $t$ of $X$ transcendental over $K$, the tree of preimages of $t$ under the iterates of $f$ carries an action of the absolute Galois group of $K(X)$: its image is the geometric iterated Galois group $G_\infty=G_\infty(K,f,t)$, a closed subgroup of the automorphism group of the tree. The fixed-point proportion $\FPP(G_\infty)$ is the limiting proportion of elements of the level-$n$ image of $G_\infty$ fixing at least one vertex at level $n$.

The level-$n$ Galois group is $A[\alpha^n]\rtimes G$, and fixing a vertex is a condition on the parameter $c\in A[\alpha^n]$ of Lemma~\ref{sec3:affine}. The paper proves one identity in two forms. The counting form expresses $\FPP(G_\infty)$ as an average over $G$ of reciprocals of fixed-point counts in $\alpha$-power torsion. The determinant form evaluates each weight as a product of $p$-adic determinants on the $\alpha$-active subspaces, and for scalar multipliers it reduces to a determinant on $H_1(A,\Z)$. From the formula we read off positivity and vanishing criteria and the T/I distinction. We prove that the fixed-point process is a martingale with an explicit limit law, and that finite groups with translations reduce to the fixed-origin case (Theorem~\ref{thm:E}). In dimension one the formula becomes explicit: ten values for scalar multipliers and nontrivial $G\le\Aut(E,0)$, for instance $5/16$ for $G=\mu_4$ with $2\mid m$. For the eight quadratic Latt\`es classes the proportions are $1/8$, $1/4$ and $5/16$, with arithmetic upper bounds over number fields (Corollary~\ref{cor:quadraticArith}). The dynamically exceptional maps other than the Latt\`es case (postcritically finite maps with a finite-type orbifold), excluded alongside it, are treated here only through the exceptional case of the Fari\~na-Asategui--Radi conjecture (Corollary~\ref{cor:conj}). All results concern the geometric group. The arithmetic analogue is a different problem, developed in a companion paper.

Fixed-point proportions go back to Odoni's work on iterated polynomial extensions \cite{Odoni}. In the arboreal setting Jones and Rouse \cite[{\S}2]{JR} formalised the invariant and computed the densities of arboreal representations in the maximal-image case \cite[Theorem~1.1]{JR}. Jones \cite{J} studied it systematically for iterated monodromy groups, with a corrigendum \cite{FJR}, and \cite{JONES} surveys the subject. Finite-index and image constraints have been obtained for cubic and unicritical polynomials \cite{BT,BDGHT} and for quadratic arboreal representations \cite{FP}. Specialisation criteria for postcritically finite rational functions appear in \cite{BGJT}, non-geometric iterated extensions in \cite{Juul}, and positive-proportion and self-similarity results in \cite{Radi2,FAR2}. In dimension one the literature covers the polynomial case \cite{FAR} and the tamely ramified non-pullback case \cite{R}. The Latt\`es case is treated here. The complex-dynamical background is surveyed in \cite{M} and treated in \cite{DH}. The higher-dimensional literature on Latt\`es maps is \cite{Din,Ro}. Arboreal representations of iterated endomorphisms of abelian algebraic groups were studied by Jones and Rouse \cite{JR}. The group-theoretic martingale criterion we use is \cite[Theorem~5.7]{BJKL}, restated in \cite[Theorem~3.5]{HZ}, and a geometric higher-dimensional criterion is \cite[Theorem~1.4]{HZ}. To our knowledge, exact fixed-point proportions for Latt\`es maps of dimension at least two have not appeared in the literature.

The counting form is the following. Let $A[\alpha^\infty]=\bigcup_{n\ge1}\ker(\alpha^n)$ be the $\alpha$-power torsion subgroup and let $\Fix_A(g)$ be the fixed-point subgroup of $g$. Then
\[
\FPP(G_\infty)=\frac1{\abs{G}}\sum_{g\in G}\frac1{\#\bigl(\Fix_A(g)\cap A[\alpha^\infty]\bigr)},
\]
with the convention that a term is zero when the intersection is infinite. The formula requires no commutativity assumption, only that $\alpha$ normalises $G$.

For an arbitrary normalising isogeny the individual weights take a local form. For each prime $p$ let $T_p(A)$ be the $p$-adic Tate module and let $W_p\subseteq T_p(A)\otimes\Q_p$ be the sum of the generalised eigenspaces of $\alpha$ on which the $p$-adic valuation of the eigenvalue is positive. We call $W_p$ the $\alpha$-active subspace. Then
\[
\#\bigl(\Fix_A(g)\cap A[\alpha^\infty]\bigr)
=\prod_p\abs{\det\nolimits_{\Q_p}(1-g\mid W_p)}_p^{-1}.
\]
For scalar $\alpha=m$ this is the closed formula of Corollary~\ref{cor:scalar}.

Three phenomena appear when $\dim A\ge2$. First, the fixed-point proportion is not a function of the abstract group $G$ alone: different embeddings of the same abstract group in $\Aut(A,0)$ can give different values, and the formula reads the multiset of weights rather than the embedding. Second, for scalar multipliers the vanishing of the average is governed by the eigenvalue-one condition. The groups concerned split into a T-class (a common fixed vector) and an I-class (no common fixed vector, while every non-identity element has eigenvalue 1). Third, positive limit values can occur outside $\{0,1,2,d\}$, where $d=\deg f$: the values $\{0,1,2\}$ already occur in earlier geometric examples, while $d$ occurs in the abstract construction of \cite[Theorem~9]{FARQ8}. This third phenomenon already has a one-dimensional precursor (Remark~\ref{rem:g1shadow}).

The one-dimensional case is treated in Section~\ref{sec:onedim}. There $A=E$ is an elliptic curve and $G\le\Aut(E,0)$ is one of the nontrivial groups $\mu_2,\mu_3,\mu_4,\mu_6$. For scalar multipliers $\alpha=m\in\Z$, Corollary~\ref{cor:values1} determines the ten values in closed form. They depend on which of $2$ and $3$ divide $m$. Corollary~\ref{cor:BJKL} shows that the limit inferior of the periodic-point proportion attains the $1/8$ bound of Bridy--Jones--Kelsey--Lodge, so the bound is sharp, and Corollary~\ref{cor:conj} settles the exceptional case of the Fari\~na-Asategui--Radi fixed-point proportion conjecture. The conjecture predicts that a tamely ramified rational function of degree $d\ge2$ whose geometric iterated Galois group is a martingale group has positive fixed-point proportion exactly when it is exceptional and not linearly conjugate to $x^{\pm d}$. The case of a non-exceptional dynamical pullback (in the sense of \cite{R}) remains open. What controls the value is the $\alpha$-adic index of $(\zeta-1)E[\alpha^\infty]$, not the complex covolume (Remark~\ref{rem:covolume}). The argument needs no fullness statement and no cardinality hypothesis on $K$ (compare Remark~\ref{sec3:rem-fullness}). For $K=\C$, the identification with the iterated monodromy group gives an alternative route to the one-dimensional values (Remark~\ref{rem:fullness-alt}).

Every group arising from an abelian-variety quotient as above is a martingale group (its fixed-point process is a martingale) in the fixed-origin case $G\le\Aut(A,0)$, by the finite-level structure and the martingale criterion of \cite{HZ,BJKL}. With the counting formula this gives the limit law of the fixed-point process: almost surely the process converges to zero or to the reciprocal of a weight. For $A=E^2$, $G=\mu_4^2$ and $m=2$ the positive limit values are $4$, $8$ and $16$, with probabilities $1/16$, $1/32$ and $1/256$.

The argument develops in four steps. Section~\ref{sec:setting} fixes the data (H) and states the main theorem in its two forms, together with the affine reduction of Theorem~\ref{thm:E}. Sections~\ref{sec:finite-level} and~\ref{sec:det} build the two tools: the identification $\Gamma_n=A[\alpha^n]\rtimes G$ of the level-$n$ Galois group (Theorem~\ref{sec3:galois}), which turns the fixed-point count into finite-group algebra, and the local index lemma, which turns the count into $p$-adic determinants. Sections~\ref{sec:zero} and~\ref{sec:process} read the formula in two directions: when the average vanishes, and how the fixed-point process behaves along the tower. Section~\ref{sec:onedim} evaluates the formula in dimension one, and Section~\ref{sec:open} lists what remains open. What is new is the determinant form of Sections~\ref{sec:setting}--\ref{sec:det}, the zero and T/I analysis, the limit law, the ten values, the attainment of $1/8$, and the exceptional case. The finite-level structure of Section~\ref{sec:finite-level} is standard (compare \cite[{\S}2]{JR} and classical covering theory), and the counting form follows from it. The arithmetic analogue over number fields is developed in a companion paper by the author (in preparation). Its unconditional values for $E:y^2=x^3-x$ over $\Q(i)$ are recorded in Section~\ref{sec:open}.

\section{Setting and main results}\label{sec:setting}

The paper rests on one standing hypothesis, (H), which fixes the field, the abelian variety, the group and the isogeny once and for all. We state the main theorem in its counting and determinant forms, and give the affine reduction that extends the formula to finite groups with translations. Two lemmas prepare the reduction and the determinant form. Their proofs can be skipped on a first reading.

\subsection{Standing hypotheses and notation}

Throughout the paper we assume the following.
\textup{(H)} $K$ is a separably closed field of characteristic $0$; $A/K$ is an abelian variety of dimension $g$; $G\le\Aut_K(A,0)$ is a finite group; $\alpha\in\End_K(A)$ is an isogeny of degree at least $2$ that normalises $G$, where $\Aut_K$ and $\End_K$ denote automorphisms and endomorphisms defined over $K$. No further commutation relation between $G$ and $\alpha$ is assumed. We write $d=\deg\alpha$ and, in the sequel, $\Aut(A,0)$ for $\Aut_K(A,0)$. Here $H_1(A,\Z)$ denotes the first singular homology of $A(\C)$ for a chosen embedding of a finitely generated subfield of definition of $(A,G,\alpha)$ into $\C$; it is a free $\Z$-module of rank $2g$ with a natural $G$-action. The module $H_1(A,\Z)$ is used below only through the determinants $\det(1-\zeta\mid H_1(A,\Z))$, which are independent of the choices. We write $\Aut(A):=\Aut_K(A)$ and $\End(A):=\End_K(A)$ for the automorphisms and endomorphisms of $A$ over $K$; every $g\in\Aut(A)$ has the form $x\mapsto L_gx+c_g$ with $L_g\in\Aut(A,0)$ and $c_g\in A(K)$, and $\pi_0(g):=L_g$ is its linear part; the map $\pi_0$ is a homomorphism.

\subsection{Abelian Latt\`es maps}

Let $A,G,\alpha$ be as in \textup{(H)}, so that $\alpha$ normalises $G$:
\begin{equation}\label{eq:normalize}
\alpha G\alpha^{-1}=G.
\end{equation}
The condition ensures that $\alpha$ descends to the quotient $X=A/G$, giving a finite map $f\colon X\to X$, the Latt\`es map associated with $(A,G,\alpha)$. Its degree is $\deg f=\deg\alpha$. For $G\le\Aut(A)$ (which may contain translations) we interpret \eqref{eq:normalize} as follows: for every $g\in G$ there is a unique $g'\in G$ with $\alpha g=g'\alpha$, and $g\mapsto g'$ is a bijection of $G$. In particular, $\alpha g\alpha^{-1}\in G$ for every $g$. Let $\pi\colon A\to X$ be the quotient map.

Let $t$ be the generic point of $X$ over $K$, so that $K(t)=K(X)$. For $n\ge1$ let $L_n$ be the set of $K(X)^{\mathrm{sep}}$-points $x$ with $f^n(x)=t$, the $d^n$ roots of the equation. The tree $\mathcal T$ with vertex set $\bigsqcup_{n\ge0}L_n$ is a rooted tree of degree $d=\deg f$, and the absolute Galois group $\Gal(K(X)^{\mathrm{sep}}/K(X))$ acts on it. Its image is the geometric iterated Galois group $G_\infty=G_\infty(K,f,t)$. Let $G_n$ be the image of $G_\infty$ in the symmetric group on $L_n$, and let $\pi_n\colon G_\infty\to G_n$ be the natural projection.

\begin{definition}[Fixed-point proportion]\label{def:fpp}
The fixed-point proportion of $G_\infty$ is
\[
\FPP(G_\infty)=\lim_{n\to\infty}\frac1{\abs{G_n}}\#\bigl\{\sigma\in G_n:\ \Fix(\sigma\mid L_n)\ne\emptyset\bigr\}.
\]
\end{definition}
The limit exists because the sequence is non-increasing; see \cite[Section 1 and (1.5)]{BJKL}. For a closed subgroup $J\le\Aut(\mathcal T)$ one defines $\FPP(J)$ as the limiting proportion of elements of the level-$n$ image of $J$ that fix at least one vertex. We call the value $\FPP(G_\infty(K,f,t))$ the geometric fixed-point proportion of $f$.

We also need finite groups $G\le\Aut(A)$ that may contain translations. When \eqref{eq:normalize} holds, $\alpha$ still descends to $A/G$ and gives a Latt\`es map. Lemma~\ref{lem:E1} below is used in the proof of Theorem~\ref{thm:E}.

\begin{lemma}[Translation part]\label{lem:E1}
Let $G\le\Aut(A)$ be finite with $\alpha G\alpha^{-1}=G$, put $T=G\cap A$, and identify $A$ with its translations. Then (i) $\alpha T=T$, (ii) $T\cap A[\alpha^\infty]=0$, and (iii) the conjugation $c_\alpha\colon G\to G$ defined by $c_\alpha(g)\alpha=\alpha g$ is an automorphism of $G$.
\end{lemma}

\begin{proof}
Conjugation by $\alpha$ induces a homomorphism $c_\alpha\colon G\to G$ with $c_\alpha(g)\alpha=\alpha g$; its image is $G$ by \eqref{eq:normalize}. Since $G$ is finite, $c_\alpha$ is bijective, hence an automorphism of $G$; this proves (iii). For $t\in T$ the translation $\tau_t\colon x\mapsto x+t$ satisfies $c_\alpha(\tau_t)=\tau_{\alpha t}$, so $\alpha T\subseteq T$. If $\alpha t=0$ for some $t\in T$, then $c_\alpha(\tau_t)=\mathrm{id}$, hence $\tau_t=\mathrm{id}$ and $t=0$, so $\alpha\mid_T$ is injective and, $T$ being finite, bijective. Thus $\alpha T=T$. Iterating, $\alpha^m\mid_T$ is bijective for every $m$, hence $T\cap\ker\alpha^m=0$, and $T\cap A[\alpha^\infty]=0$.
\end{proof}

\subsection{The \texorpdfstring{$\alpha$}{alpha}-adic Tate module and active subspaces}

The $\alpha$-adic Tate module is $T_\alpha=T_\alpha(A)=\varprojlim_{n}(A[\alpha^n],\times\alpha)$. For a prime $p$, let $T_p(A)$ be the $p$-adic Tate module, put $V_p=T_p(A)\otimes_{\Z_p}\Q_p$ and $M_p=\alpha\mid V_p$, and let $V_p=\bigoplus_{\lambda}V_{p,\lambda}$ be the primary decomposition over the irreducible factors $\lambda$ of the characteristic polynomial of $M_p$ over $\Q_p$. All roots of a fixed $\lambda$ have the same $p$-adic valuation (Lemma~\ref{sec4:val}). We put
\[
W_p=\bigoplus_{v_p(\lambda)>0}V_{p,\lambda},
\qquad
U_p=\bigoplus_{v_p(\lambda)=0}V_{p,\lambda}.
\]
We call $W_p$ the $\alpha$-active subspace at $p$ and $U_p$ the inactive subspace. Then $V_p=W_p\oplus U_p$, and $W_p=0$ if and only if $p\nmid\deg\alpha$ (Lemma~\ref{sec4:act}). We write $W_p=W_p(\alpha)$ when the isogeny is to be named, and $|\cdot|_p$ for the $p$-adic absolute value normalised by $|p|_p=1/p$.

\begin{lemma}[Centralising power]\label{sec2:beta}
Assume \textup{(H)} and put $k=\abs{\Aut(G)}$. Then $\beta:=\alpha^k$ is an isogeny of degree $(\deg\alpha)^k$ satisfying \eqref{eq:normalize}, and (a) $\beta$ centralises $G$, (b) $A[\beta^\infty]=A[\alpha^\infty]$, (c) $W_p(\beta)=W_p(\alpha)$ for every prime $p$, where $W_p(\beta)$ is the $\beta$-active subspace.
\end{lemma}

\begin{proof}
By Lemma~\ref{lem:E1}(iii), the conjugation $c_\alpha\colon G\to G$ is an automorphism of $G$, so its order divides $k=\abs{\Aut(G)}$. Since $c_{\alpha^m}=c_\alpha^m$ for every $m$, we have $c_\beta=c_\alpha^k=\mathrm{id}$; hence $\beta$ centralises $G$, and $\beta G\beta^{-1}=G$ because $\alpha^mG\alpha^{-m}=G$ for every $m$. For (b), $\ker(\beta^n)=\ker(\alpha^{kn})$ and the multiples of $k$ are cofinal, so $A[\beta^\infty]=A[\alpha^\infty]$. For (c) let $V_{p,\lambda}$ be a generalised eigenspace of $M_p=\alpha\mid V_p$ for the factor $\lambda$. The eigenvalues of $M_p^k$ on it have valuation $k\,v_p(\lambda)$, which is positive if and only if $v_p(\lambda)>0$. If two factors $\lambda,\lambda'$ contribute to the same generalised eigenspace of $M_p^k$, then some root $\alpha$ of $\lambda$ and some root $\beta$ of $\lambda'$ satisfy $\alpha^k=\beta^k$; since all roots of a fixed factor have the same valuation (Lemma~\ref{sec4:val}), $k\,v_p(\lambda)=v_p(\alpha^k)=v_p(\beta^k)=k\,v_p(\lambda')$. Hence every generalised eigenspace of $M_p^k$ is entirely positive or entirely non-positive, and the generalised eigenspaces with positive valuation are exactly the sums of the $V_{p,\lambda}$ with $v_p(\lambda)>0$, so $W_p(\beta)=W_p(\alpha)$.
\end{proof}

\subsection{Main theorem}

Lemma~\ref{sec2:beta} shows that the centrality used in Section~\ref{sec:det} can always be arranged. Parts (i) and (ii) are one statement read twice: first as a count in $\alpha$-power torsion, then as a product of $p$-adic determinants on those subspaces. Sections~\ref{sec:finite-level} and~\ref{sec:det} prove the two parts.

\begin{theorem}[Main theorem]\label{thm:main}
Let $A,G,\alpha$ be as in \textup{(H)}.
\begin{enumerate}
\item[(i)] (Counting form) The fixed-point proportion is
\[
\FPP(G_\infty)=\frac1{\abs{G}}\sum_{g\in G}\frac1{\#\bigl(\Fix_A(g)\cap A[\alpha^\infty]\bigr)},
\]
where a term is zero when the intersection is infinite.
\item[(ii)] (Determinant form) For every $g\in G$,
\[
\#\bigl(\Fix_A(g)\cap A[\alpha^\infty]\bigr)
=\prod_p\abs{\det\nolimits_{\Q_p}(1-g\mid W_p)}_p^{-1},
\]
where a factor is infinite when the determinant vanishes. Consequently
\[
\FPP(G_\infty)=\frac1{\abs{G}}\sum_{g\in G}\prod_p\abs{\det\nolimits_{\Q_p}(1-g\mid W_p)}_p.
\]
\end{enumerate}
\end{theorem}

For $g\in G$ put $w_\alpha(g):=1/\#(\Fix_A(g)\cap A[\alpha^\infty])\in[0,1]$, with $w_\alpha(g)=0$ when the intersection is infinite. Then $\FPP=(1/|G|)\sum_{g\in G}w_\alpha(g)$, and by (ii) $w_\alpha(\zeta)=\kappa(\zeta)^{-1}$ with $\kappa(\zeta):=\prod_p\abs{\det(1-\zeta\mid W_p(\alpha))}_p^{-1}\in\N\cup\{\infty\}$. Both forms hold under \textup{(H)} alone, and the determinant form reduces to the central case in Section~\ref{sec:det} by Lemma~\ref{sec2:beta}. The identities $w_\alpha=\kappa^{-1}=k_\infty^{-1}$ are used below, with $k_\infty$ as in Section~\ref{sec:finite-level}; in the scalar case $w_\alpha(\zeta)=\rho_m(\zeta)=\abs{D(\zeta)}_m$, with $\rho_m$ and $D$ as in Section~\ref{sec:process} and Theorem~\ref{thm:zero}, and $\abs{\cdot}_m$ as in Section~\ref{sec:onedim}.

\begin{corollary}[Scalar case]\label{cor:scalar}
Assume $\alpha=m\in\Z$ with $\abs{m}\ge2$. Then $W_p=V_p$ for $p\mid m$ and $W_p=0$ otherwise, and
\[
\FPP(G_\infty)=\frac1{\abs{G}}\sum_{1\ne\zeta\in G}
\prod_{p\mid m}\abs{\det(1-\zeta\mid H_1(A,\Z))}_p .
\]
\end{corollary}

\begin{proof}
For scalar $\alpha=m$ the isogeny is central. The eigenvalues of $\alpha$ on $V_p$ equal $m$, so $W_p=V_p$ for $p\mid m$ and $W_p=0$ otherwise by Lemma~\ref{sec4:act}(iii). Theorem~\ref{thm:main}(ii) gives $w_\alpha(\zeta)=\prod_{p\mid m}\abs{\det(1-\zeta\mid V_p)}_p$. Since $\zeta$ has finite order, $\det(1-\zeta\mid V_p)=\det(1-\zeta\mid H_1(A,\Z))$, and inserting this into Theorem~\ref{thm:main}(i), with $\zeta=1$ omitted, gives the formula.
\end{proof}

\begin{corollary}[Positivity, upper bound, invariance]\label{cor:pos}
\begin{enumerate}
\item[(a)] $\FPP(G_\infty)>0$ if and only if there exists $g\ne1$ in $G$ with $\#(\Fix_A(g)\cap A[\alpha^\infty])<\infty$.
\item[(b)] $\FPP(G_\infty)\le(\abs{G}-1)/\abs{G}$.
\item[(c)] The value $\FPP(G_\infty)$ depends only on the $\Z[G]$-module $H_1(A,\Z)$ and the action of $\alpha$; in the scalar case it depends only on the $\Z[G]$-module $H_1(A,\Z)$ and on $m$.
\end{enumerate}
The bound in (b) is attained: for $G=\mu_4$ acting on an elliptic curve and $m$ odd, Corollary~\ref{cor:scalar} gives $\FPP=3/4=(\abs{G}-1)/\abs{G}$.
\end{corollary}

\begin{proof}
(a) By Theorem~\ref{thm:main}(i), $\FPP$ is the average of the non-negative weights $w_\alpha(g)$, so it vanishes exactly when $\Fix_A(g)\cap A[\alpha^\infty]$ is infinite for every $g\ne1$.

(b) The identity has infinite fixed-point set in $A[\alpha^\infty]$, so $w_\alpha(1)=0$, while $0\le w_\alpha(g)\le1$, so $\FPP\le(\abs{G}-1)/\abs{G}$, with equality for $G=\mu_4$ and $m$ odd by Corollary~\ref{cor:scalar}.

(c) The sets $\Fix_A(g)\cap A[\alpha^\infty]$ are determined by the $\Z[G]$-module $H_1(A,\Z)$ and $\alpha$: there is a $\Z$-module isomorphism $A[\alpha^n]\cong\alpha^{-n}H_1(A,\Z)/H_1(A,\Z)$, with $\alpha^{-n}H_1(A,\Z)$ the preimage of $H_1(A,\Z)$ in $H_1(A,\Z)\otimes_\Z\Q$ under $\alpha^n$, and the action of $g$ is transported through the conjugate endomorphism $\alpha^n g\alpha^{-n}$. In the scalar case Corollary~\ref{cor:scalar} leaves only $H_1(A,\Z)$ and $m$.
\end{proof}

\begin{example}\label{ex:sameorder}
Let $A=E^2$ and $m=2$. For $G=\langle\mathrm{swap}\rangle$ the two factors are exchanged, so every non-identity element has zero weight and $\FPP=0$. For $G=\{\pm I\}$ the same ambient variety and the same multiplier give $\FPP=1/32$. Both groups have order $2$, so the fixed-point proportion is not a function of the isomorphism class of $G$ alone. Different embeddings of one abstract group can give different values: for $E:y^2=x^3-x$ (so $\mu_4\le\Aut(E)$), the six order-four cyclic subgroups of $\mu_4\times\mu_4\le\Aut(E^2)$ give exactly three values, namely $0$, $1/16$ and $9/64$. The values are computed from Corollary~\ref{cor:scalar}.
\end{example}

\begin{remark}\label{rem:noncentral}
Neither form needs centrality. The counting form is proved in Section~\ref{sec:finite-level}. For the determinant form Lemma~\ref{sec2:beta} replaces $\alpha$ by $\beta=\alpha^k$, which centralises $G$ and leaves $A[\alpha^\infty]$ and every $W_p$ unchanged, after which Section~\ref{sec:det} applies. For an isogeny that only normalises $G$, the linear part at level $n$ is $\gamma_n=\alpha^{-n}\zeta\alpha^n$ rather than $\zeta$ (Appendix~\ref{app:twist}), while the weight is still the product above.
\end{remark}

The main theorem assumes $G\le\Aut_K(A,0)$. Theorem~\ref{thm:E} below, whose proof is given in Section~\ref{sec3:affine-general}, reduces the fixed-point proportion for a finite group that may contain translations to the fixed-origin case. We write $\FPP(A/G,\alpha)$ for the fixed-point proportion of the associated Latt\`es map. This value is intrinsic to $(A,G,\alpha)$: with $t$ the generic point, $G_\infty(K,f,t)$ is determined by $(K,f)$.

\begin{theorem}[Affine automorphism groups]\label{thm:E}
Assume \textup{(H)}, except that $G\le\Aut_K(A)$ is finite and may contain translations. Put $H_0=\pi_0(G)$. Then $\alpha H_0\alpha^{-1}=H_0$. Moreover, with the convention $1/\infty=0$,
\[
\FPP(A/G,\alpha)=\frac1{\abs{H_0}}\sum_{h\in H_0}\frac1{\#\bigl(\Fix_A(h)\cap A[\alpha^\infty]\bigr)}
=\frac1{\abs{H_0}}\sum_{h\in H_0}\prod_p\abs{\det(1-h\mid W_p(\alpha))}_p,
\]
and this common value equals $\FPP(A/H_0,\alpha)$.
\end{theorem}

\begin{remark}
When $\zeta\in G$ has a nonzero translation part, its action on $H_1(A,\Z)$ and on $V_p$ is induced by its linear part $L_\zeta:=\pi_0(\zeta)$. Its fixed-point set, however, is not determined by that linear part: a nonzero translation has no fixed points, while the identity element fixes all of $A$. Theorem~\ref{thm:E} therefore computes the weights from the linear parts, and $\FPP$ depends only on $H_0=\pi_0(G)$.
\end{remark}

\subsection{Hypotheses at a glance}

The table below lists the additional hypotheses, beyond (H), used by the principal results of Sections~\ref{sec:setting}--\ref{sec:process}. The results of Section~\ref{sec:onedim} state their hypotheses individually.

\begin{center}
\small
\begin{tabular}{@{}ll@{}}
\toprule
Result & Additional hypothesis \\
\midrule
Theorem \ref{thm:main}(i) & none beyond (H) \\
Theorem \ref{thm:main}(ii) & none beyond (H) \\
Corollary \ref{cor:scalar} & $\alpha=m\in\Z$ \\
Corollary \ref{cor:pos} & none beyond (H) \\
Theorem \ref{thm:E} & no fixed-origin assumption on $G$; $\alpha$ normalises $G$ \\
Theorem \ref{sec3:galois} & none beyond (H) \\
Theorem \ref{thm:zero} & (i)--(iii): none beyond (H); (iv) $\alpha=m\in\Z$ \\
Proposition \ref{prop:dim2} & $\dim A=2$, $\alpha=m\in\Z$ \\
Theorem \ref{thm:TI} & $\alpha=m\in\Z$ \\
Theorem \ref{thm:Ifamilies} & (a)--(c) $\alpha=m\in\Z$; (d) arbitrary normalising isogeny \\
Theorem \ref{thm:martingale} & none beyond (H) \\
Theorem \ref{thm:limitlaw} & none beyond (H) \\
Corollary \ref{cor:Q8} & $\alpha=m\in\Z$ \\
\bottomrule
\end{tabular}
\end{center}

\section{Finite-level structure}\label{sec:finite-level}

The preimage tree of $f$ is governed by the $\alpha$-power torsion of $A$, and this is what turns the fixed-point count into finite-group algebra. Only the hypothesis \textup{(H)} is used here; in particular $G$ need not commute with $\alpha$. Section~\ref{sec3:affine-general} records the small modification needed when $G$ contains translations. The counting form of the main theorem, Theorem~\ref{thm:main}(i), is proved here. Throughout, \(F=K(X)\) is the function field of \(X\), \(t\) its generic point, and \(P\) the generic point of \(A\) over \(K\), so that \(K(P)=K(A)\) and \(\pi(P)=t\). Choose $y_0=P$ and then, recursively, points $y_n\in A(F^{\mathrm{sep}})$ with $\alpha y_n=y_{n-1}$ for $n\ge1$; then $\alpha^n y_n=P$. The base field is \(F\), and nothing in the argument depends on $\dim A$.

\subsection{The covering map and its deck group}
\label{sec3:deck}

We begin with the covering whose deck group acts on level $n$. Put \(\psi_n:=\pi\circ\alpha^n\colon A\to X\), so that
\(\psi_n=f^n\circ\pi\). Since \(\deg\pi=\abs{G}\) and \(\deg\alpha^n=d^n\),
the morphism \(\psi_n\) is finite and separable of degree \(d^n\abs{G}\).

\begin{theorem}[Galois structure]\label{sec3:galois}
The field \(K(y_n)\) is a Galois extension of \(F\); over \(K(A)\) its
Galois group is \(A[\alpha^n]\), and over \(F\) its Galois group is
\[
\Gamma_n:=A[\alpha^n]\rtimes G ,
\]
where \(G\) acts on \(A[\alpha^n]\) through its natural action. Moreover
\(\Gamma_n\) is canonically isomorphic to the deck group
\[
\Deck(\psi_n):=\{\tau\in\Aut(A):\ \psi_n\circ\tau=\psi_n\}
 =\bigl\{x\mapsto gx+h_0\ \bigm|\ g\in G,\ h_0\in A[\alpha^n]\bigr\},
\]
and \(\abs{\Gamma_n}=d^n\abs{G}=\deg\psi_n\); the subfield \(K(A)=K(P)\) is the fixed field of \(A[\alpha^n]\) in \(K(y_n)\). Here the group law is
\((h_0,g)(h_0',g')=(h_0+g h_0',gg')\).
\end{theorem}

\begin{proof}
First \([K(y_n):K(A)]=d^n\), because \(y_n\) is a root of the equation
\(\alpha^n y=P\) over the generic point \(P\) and the isogeny
\(\alpha^n\colon A\to A\) is a Galois covering with group \(A[\alpha^n]\)
(characteristic \(0\)). Since \(K(y_n)\supseteq K(P)=K(A)\) and
\([K(A):F]=\abs{G}\), we get \([K(y_n):F]=d^n\abs{G}\).

Let $\tau\in\Aut(A)$ satisfy $\psi_n\circ\tau=\psi_n$. Every automorphism
of \(A\) has the form \(\varphi(x)=gx+h_0\) with \(g\in\Aut(A,0)\) and
\(h_0\in A(K)\); in particular, \(\tau\) has this form. The condition reads
\[
\pi\bigl(\alpha^n(gx+h_0)\bigr)=\pi(\alpha^n x)\qquad\text{for all }x\in A.
\]
Write \(g_n:=\alpha^n g\alpha^{-n}\in\Aut(A,0)\) and put
\(y:=\alpha^n x\). Since \(\alpha^n\) is surjective, the condition is
\[
\pi(g_n y+\alpha^n h_0)=\pi(y)\qquad\text{for all }y\in A .
\]
Taking \(y=0\) gives \(\pi(\alpha^nh_0)=\pi(0)\); since the fibre of
\(\pi\) over \(\pi(0)\) is the orbit \(G\cdot0=\{0\}\), we get
\(\alpha^nh_0=0\). Then \(\pi(g_ny)=\pi(y)\) for all \(y\), which means
\(g_n\in G\) (as \(G\) is the group of automorphisms of \(A\) over \(X\)). Hence
\(g\in G\), because \(\alpha^nG\alpha^{-n}=G\), and
\(h_0\in A[\alpha^n]\). Conversely,
if \(g\in G\) and \(h_0\in A[\alpha^n]\), then
\(\psi_n(gx+h_0)=\pi(g_n\alpha^nx+\alpha^nh_0)=\pi(g_n\alpha^nx)
=\pi(\alpha^nx)=\psi_n(x)\). This proves the description of
\(\Deck(\psi_n)\); its order is \(d^n\abs{G}\).

Finally, the Galois group. By the preceding paragraph,
\(\Deck(\psi_n)=\{x\mapsto gx+h_0:\ g\in G,\ h_0\in A[\alpha^n]\}\) has order
\(d^n\abs G=\deg\psi_n\). It acts transitively on the fibre of \(\psi_n\) over the
generic point \(t\): every point \(x\) of that fibre satisfies \(\alpha^nx=\gamma P\)
for some \(\gamma\in G\), and choosing \(\gamma_1\in G\) with
\(\alpha^n\gamma_1=\gamma\alpha^n\) gives \(x=\gamma_1(y_n)+h\) with
\(h\in A[\alpha^n]\) (all points of \(A[\alpha^n]\) are \(K\)-rational, since \(K\)
is separably closed); for two such points \(\gamma_1(y_n)+h\) and
\(\gamma_1'(y_n)+h'\), the element
\(x\mapsto\gamma_1'\gamma_1^{-1}x+(h'-\gamma_1'\gamma_1^{-1}h)\) of
\(\Deck(\psi_n)\) sends the first to the second, because \(A[\alpha^n]\) is
\(G\)-stable. Since \(\psi_n\) is finite and separable (characteristic \(0\)) and
its deck group has order \(d^n\abs G=\deg\psi_n\) and acts transitively on the
generic fibre, \(\psi_n\) is a finite separable Galois covering with deck group
\(\Deck(\psi_n)\) (standard theory of finite separable coverings in
characteristic \(0\)). The function field of the generic fibre of \(\psi_n\) over
\(F=K(X)\) is \(K(y_n)\) when \(\psi_n(y_n)=t\): the function field of the generic
fibre of a finite morphism of integral varieties is the function field of its
total space, and the point \(y_n\) identifies it with \(K(y_n)\). Hence
\(K(y_n)/F\) is a Galois extension with Galois group
\(\Deck(\psi_n)\cong A[\alpha^n]\rtimes G=\Gamma_n\), the parametrisation and
the group law being those of the preceding paragraph. Under this identification
the elements with \(g=1\) send \(y_n\) to \(y_n+h_0\), hence fix
\(P=\alpha^ny_n\) and therefore \(K(A)=K(P)\); since
\([K(y_n):K(A)]=d^n=\abs{A[\alpha^n]}\), the fixed field of \(A[\alpha^n]\) is
\(K(A)\).
\end{proof}

\subsection{Level vertices and the affine action}
\label{sec3:action}

Using the group, we label the level-$n$ vertices by the torsion group $A[\alpha^n]$. The next two lemmas establish the labelling and the action that follows from it.

\begin{lemma}[Level identification]\label{sec3:bijection}
The map
\[
\Phi_n\colon A[\alpha^n]\longrightarrow f^{-n}(t),\qquad
h\longmapsto \pi(y_n+h)
\]
is a bijection, and under \(\Phi_n\) the parent map \(f\) corresponds to
\(h\mapsto\alpha h\).
\end{lemma}

\begin{proof}
Suppose that $\pi(y_n+h)=\pi(y_n+h')$ for some $h,h'\in A[\alpha^n]$. Then there is $\gamma\in G$ with
\[
y_n+h'=\gamma(y_n+h).
\]
Rearranging gives $h'-\gamma h=(\gamma-1)y_n$. The left-hand side is a torsion point, hence is defined over $K$. Since $K(y_n)$ is finite over $F=K(X)$, the point $y_n$ is generic over $K$. If $\gamma\ne1$, then $\gamma-1$ is a nonzero endomorphism, its image $B$ is a positive-dimensional abelian subvariety of $A$, and $(\gamma-1)y_n$ is the generic point of $B$. A torsion point is not the generic point of a positive-dimensional variety, so $\gamma=1$ and then $h=h'$. Hence $\Phi_n$ is injective.

For surjectivity, let $z\in f^{-n}(t)$ and choose \(x\in A(F^{\mathrm{sep}})\) with $\pi(x)=z$ (possible since $\pi$ is finite and separable in characteristic \(0\)). Then $\alpha^nx\in\pi^{-1}(t)=G\cdot P$, so $\alpha^nx=\gamma P$ for some $\gamma\in G$. Put $\gamma':=\alpha^{-n}\gamma^{-1}\alpha^n\in G$. Then
\[
\alpha^n(\gamma'x)=\gamma^{-1}\alpha^nx=P,
\]
so $\gamma'x=y_n+h$ for some $h\in A[\alpha^n]$. Thus
\[
z=\pi(x)=\pi(\gamma'x)=\pi(y_n+h)=\Phi_n(h),
\]
where the second equality uses $\pi\circ\gamma'=\pi$ and the last one is the definition of $\Phi_n$. Hence $\Phi_n$ is a bijection. The parent formula follows from $f(\pi(y_n+h))=\pi(\alpha y_n+\alpha h)=\pi(y_{n-1}+\alpha h)$.
\end{proof}

\begin{lemma}[Absolute action and fixed vertices]\label{sec3:affine}
Let $\sigma\in\Gal(K(y_n)/F)$. Let $\zeta=\sigma\mid K(A)\in G$ and put
\[
\gamma_n:=\alpha^{-n}\zeta\alpha^n\in G .
\]
Then there is a unique $c\in A[\alpha^n]$ such that
\[
\sigma(y_n)=\gamma_n y_n+c .
\]
The map $\sigma\mapsto(\zeta,c)$ is a bijection
$\Gal(K(y_n)/F)\to G\times A[\alpha^n]$. In the coordinates
$v_h=\Phi_n(h)=\pi(y_n+h)$ of Lemma~\ref{sec3:bijection}, the induced
permutation of $f^{-n}(t)$ is
\[
h\longmapsto \gamma_n^{-1}(h+c).
\]
Moreover $\sigma$ fixes at least one vertex of level $n$ if and only if
$c\in(\gamma_n-1)A[\alpha^n]$; in that case it fixes exactly
$\#\bigl(A[\alpha^n]\cap\Fix_A(\gamma_n)\bigr)$ vertices.
\end{lemma}

\begin{proof}
Since $K(A)/F$ is Galois with group $G$, the automorphism $\sigma$ acts on $A(F^{\mathrm{sep}})$ through $\zeta$, so $\sigma(P)=\zeta P$. With $P=\alpha^n y_n$, applying $\alpha^n$ to
$\sigma(y_n)$ gives
\[
\alpha^n\sigma(y_n)=\zeta P=\zeta\alpha^n y_n=\alpha^n\gamma_n y_n,
\]
because $\alpha^n\gamma_n=\zeta\alpha^n$. Hence
$c:=\sigma(y_n)-\gamma_n y_n\in\ker\alpha^n=A[\alpha^n]$.

For the bijection, the map is injective: $\sigma$ is determined by its
restriction to $K(A)$ and by its value at $y_n$, and $K(y_n)$ is generated
by $K(A)$ and $y_n$. Both source and target have $d^n|G|$ elements by
Theorem~\ref{sec3:galois}, so injectivity implies bijectivity.

For a vertex $v_h=\Phi_n(h)$ we have, since $\sigma$ fixes $K$ pointwise
and $h$ is torsion,
\[
\sigma(v_h)=\pi\bigl(\sigma(y_n)+h\bigr)
=\pi\bigl(\gamma_n y_n+c+h\bigr)
=\pi\bigl(y_n+\gamma_n^{-1}(c+h)\bigr)
=\Phi_n\bigl(\gamma_n^{-1}(c+h)\bigr),
\]
where the penultimate equality uses $\pi\circ\gamma_n=\pi$. Thus the
permutation is as displayed. It fixes $v_h$ exactly when
$h=\gamma_n^{-1}(h+c)$, that is, $(\gamma_n-1)h=c$. If this holds, the
fixed vertices are the solutions of $(\gamma_n-1)h=c$ in $A[\alpha^n]$,
and their number is $\#\ker(\gamma_n-1\mid A[\alpha^n])
=\#(A[\alpha^n]\cap\Fix_A(\gamma_n))$.
\end{proof}

\subsection{The finite-level proportion}
\label{sec3:fppn}

The labelling turns the level-$n$ fixed-point count into a count in the finite group $A[\alpha^n]\rtimes G$. For \(n\ge1\) let \(\FPP_n\) be the proportion of elements of
\(\pi_n(G_\infty(K,f,t))\) that fix at least one vertex at level
\(n\).

\begin{proposition}[Finite-level formula]\label{sec3:fppn-formula}
For every \(n\ge1\),
\[
\FPP_n
 =\frac{1}{\abs{G}}\sum_{g\in G}
   \frac{1}{\#\bigl(A[\alpha^n]\cap\Fix_A(g)\bigr)}
 =\frac{1}{d^n\abs{G}}\sum_{g\in G}\#\bigl((g-1)A[\alpha^n]\bigr).
\]
\end{proposition}

\begin{proof}
By Theorem~\ref{sec3:galois} the splitting field \(K_n\) of
\(f^n(x)-t\) is the fixed field in \(K(y_n)\) of the kernel
\[
H_n:=\{\sigma\in\Gamma_n:\ \sigma \text{ acts trivially on } f^{-n}(t)\},
\]
so \(G_n=\Gal(K_n/F)\cong\Gamma_n/H_n\). Since \(H_n\) acts trivially on
the level-\(n\) vertices, the proportion in \(G_n\) equals the proportion
in \(\Gamma_n\). By Lemma~\ref{sec3:affine}, for fixed $\zeta\in G$ the elements
\(\sigma\) with $\sigma\mid K(A)=\zeta$ and with a fixed vertex are
those whose parameter $c$ lies in $(\gamma_n-1)A[\alpha^n]$, where
\(\gamma_n=\alpha^{-n}\zeta\alpha^n\); there are
\(\#((\gamma_n-1)A[\alpha^n])\) such elements. Since
$\zeta\mapsto\gamma_n$ is a bijection of $G$, summing over $\zeta$ gives
\(\sum_{g\in G}\#((g-1)A[\alpha^n])\), while
\(\abs{\Gamma_n}=d^n\abs{G}\). This gives the second expression. For the
first, note that \(g-1\) is an endomorphism of the finite group
\(A[\alpha^n]\), so
\(\#((g-1)A[\alpha^n])=d^n/\#(A[\alpha^n]\cap\Fix_A(g))\).
\end{proof}

\subsection{The limit and the counting form}
\label{sec3:limit}

Letting $n\to\infty$ turns the finite-level formula into the counting form of Theorem~\ref{thm:main}(i).

\begin{lemma}[Monotone limit]\label{sec3:monotone}
For \(g\in G\) put \(k_n(g):=\#\bigl(A[\alpha^n]\cap\Fix_A(g)\bigr)\) and
\(k_\infty(g):=\#\bigl(\Fix_A(g)\cap A[\alpha^\infty]\bigr)\). Then
\(k_1(g)\le k_2(g)\le\cdots\) and
\(\lim_n k_n(g)=k_\infty(g)\in\N\cup\{\infty\}\). Consequently
\(\lim_n\FPP_n=\frac1{\abs{G}}\sum_{g\in G}1/k_\infty(g)\).
\end{lemma}

\begin{proof}
Since \(A[\alpha^n]=\ker(\alpha^n)\subseteq\ker(\alpha^{n+1})
=A[\alpha^{n+1}]\), the sets \(A[\alpha^n]\cap\Fix_A(g)\) form an
increasing chain whose union is \(\Fix_A(g)\cap A[\alpha^\infty]\). The
first assertion follows; the second follows from
Proposition~\ref{sec3:fppn-formula} and finiteness of \(G\).
\end{proof}

\begin{proof}[Proof of Theorem~\ref{thm:main}(i)]
Let \(E_n\subseteq G_\infty\) be the set of elements fixing at least one
vertex at level \(n\). Since the parent of a fixed vertex is fixed, the
sets \(E_n\) decrease, and their intersection is the set of elements fixing a
vertex at every level. By continuity of the Haar measure,
\(\FPP=\lim_n\mu(E_n)=\lim_n\FPP_n\). Combining this with
Lemma~\ref{sec3:monotone} gives the counting form.
\end{proof}

\begin{remark}\label{sec3:rem-fullness}
Neither the identification of \(G_\infty\) with the iterated monodromy group nor a fullness statement for the arboreal image is used. The proof uses only the Galois group of \(K(y_n)/F\). In particular the conclusion holds for every separably closed \(K\) of characteristic \(0\), without a cardinality assumption.
\end{remark}

\subsection{The affine case}
\label{sec3:affine-general}

We now allow $G\le\Aut(A)$ to contain translations. Write $g(x)=L_gx+c_g$ with $L_g=\pi_0(g)$ and $c_g=g(0)$, and put $T=G\cap A$ and $H_0=\pi_0(G)\le\Aut(A,0)$.

\begin{lemma}[Level structure for affine groups]\label{sec3:affine-level}
Assume \textup{(H)}, except that $G\le\Aut(A)$ is finite and may contain translations. For $n\ge1$ and $L\in\Aut(A,0)$ put $k_n(L):=\#\bigl(A[\alpha^n]\cap\Fix_A(L)\bigr)$. Then $\FPP_n=(1/\abs G)\sum_{g\in G}1/k_n(L_g)$.
\end{lemma}

\begin{proof}
The proportion in $G_n$ equals that in $\Gamma_n:=\Gal(K(y_n)/F)$, as in the proof of Proposition~\ref{sec3:fppn-formula}. First, as in the first paragraph of the proof of Theorem~\ref{sec3:galois}, $[K(y_n):F]=d^n\abs G$; the fibre of $\psi_n=\pi\circ\alpha^n$ over $t$ is $\bigcup_{\gamma\in G}\{x:\alpha^nx=\gamma P\}$, and for $\gamma_1\in G$ with $\alpha^n\gamma_1=\gamma\alpha^n$ (Lemma~\ref{lem:E1}(iii)) the corresponding block is $\gamma_1(y_n)+A[\alpha^n]$. Since $\gamma_1$ is defined over $K$ and $A[\alpha^n]\subset A(K)$, each point of this block lies in $A(K(y_n))$. The fibre has $d^n\abs G=[K(y_n):F]$ points and consists of conjugates of $y_n$, so $K(y_n)/F$ is Galois.

Second, as in Lemma~\ref{sec3:affine}, $\sigma\mapsto(\zeta,c)$ is a bijection $\Gamma_n\to G\times A[\alpha^n]$, where $\zeta=\sigma\mid_{K(A)}$ and $\sigma(y_n)=\gamma_ny_n+c$ with $\alpha^n\gamma_n=\zeta\alpha^n$; here $c\in A[\alpha^n]$ because $\alpha^n\sigma(y_n)=\zeta\alpha^ny_n=\alpha^n\gamma_ny_n$.

Third, $\Phi_n\colon A[\alpha^n]\to f^{-n}(t)$, $h\mapsto\pi(y_n+h)$, is a bijection. If $\pi(y_n+h)=\pi(y_n+h')$, then $y_n+h'=\gamma(y_n+h)$ for some $\gamma\in G$: for $L_\gamma=1$ we have $c_\gamma\in T$ and $h'-h=c_\gamma\in A[\alpha^n]\cap T=0$ by Lemma~\ref{lem:E1}(ii), while for $L_\gamma\ne1$ the generic point $\gamma(y_n)-y_n$ of the positive-dimensional translate $(L_\gamma-1)A+c_\gamma$ would be the torsion point $h'-L_\gamma h$, impossible. Surjectivity follows as in Lemma~\ref{sec3:bijection}, with $\gamma'$ chosen so that $\alpha^n\gamma'=\gamma^{-1}\alpha^n$.

Finally, as in Lemma~\ref{sec3:affine}, $\sigma$ acts by $h\mapsto L_{\gamma_n}^{-1}(c+h)$, since $\sigma(\Phi_n(h))=\pi(\gamma_ny_n+c+h)=\Phi_n(L_{\gamma_n}^{-1}(c+h))$, where $\pi\circ\gamma_n=\pi$, $\gamma_n(y_n+u)=\gamma_n(y_n)+L_{\gamma_n}u$ for every $u\in A$, and $L_{\gamma_n}^{-1}(c+h)\in A[\alpha^n]$ because $\alpha^nL_{\gamma_n}=L_\zeta\alpha^n$. Thus $\sigma$ fixes a vertex exactly when $c\in(L_{\gamma_n}-1)A[\alpha^n]$, and then it fixes $k_n(L_{\gamma_n})$ vertices. Since $\#\bigl((L-1)A[\alpha^n]\bigr)=d^n/k_n(L)$ for the endomorphism $L-1$ of the finite group $A[\alpha^n]$, this gives
\[
\FPP_n=\frac1{d^n\abs G}\sum_{\zeta\in G}\#\bigl((L_{\gamma_n}-1)A[\alpha^n]\bigr)
=\frac1{\abs G}\sum_{\zeta\in G}\frac1{k_n(L_{\gamma_n})}
=\frac1{\abs G}\sum_{g\in G}\frac1{k_n(L_g)},
\]
where the last equality is the relabelling $g=\gamma_n$.
\end{proof}

\begin{proof}[Proof of Theorem~\ref{thm:E}]
By Lemma~\ref{sec3:affine-level} and the monotonicity argument in the proof of Theorem~\ref{thm:main}(i), $\FPP=\lim_n\FPP_n$, and $k_n(L)$ increases to $k_\infty(L):=\#\bigl(\Fix_A(L)\cap A[\alpha^\infty]\bigr)$ for $L\in\Aut(A,0)$. Hence
\[
\begin{aligned}
\FPP&=\frac1{\abs G}\sum_{g\in G}\frac1{k_\infty(L_g)}\\
&=\frac1{\abs{H_0}}\sum_{h\in H_0}\frac1{k_\infty(h)}\\
&=\frac1{\abs{H_0}}\sum_{h\in H_0}\prod_p\abs{\det(1-h\mid W_p(\alpha))}_p\\
&=\FPP(A/H_0,\alpha),
\end{aligned}
\]
where the second equality uses that $g\mapsto L_g$ maps $G$ onto $H_0$ with fibres of size $\abs T$, and the last two use Theorem~\ref{thm:main} for $(A,H_0,\alpha)$, which is legitimate since $H_0\le\Aut(A,0)$ is finite and $\alpha H_0\alpha^{-1}=H_0$ (take linear parts in $\alpha G\alpha^{-1}=G$: $\alpha L_g\alpha^{-1}=L_{\alpha g\alpha^{-1}}\in H_0$).
\end{proof}

\section{The determinant formula}\label{sec:det}

The counting form leaves the weights $w_\alpha(g)$ as reciprocal fixed-point counts. The question is what those weights are, and the answer is local: each weight is a product over primes of $p$-adic determinants on the $\alpha$-active subspace, which is the determinant form of Theorem~\ref{thm:main}(ii).

\subsection{Active subspaces}
\label{sec4:active}

We first record how $\alpha$ acts on each $p$-adic Tate module: the valuations of its eigenvalues are well defined, and they single out the active part. Let \(p\) be a prime. Put $V_p=T_p(A)\otimes_{\Z_p}\Q_p$, which is
$2\dim A$-dimensional over \(\Q_p\), and
\(M_p=\alpha\otimes\Q_p\in\GL(V_p)\).

\begin{lemma}\label{sec4:val}
Let \(\chi_p(t)=\det(t-M_p)\in\Z_p[t]\), and let
\(\lambda\in\Q_p[t]\) be a monic irreducible factor of \(\chi_p\). Then all
roots of \(\lambda\) in \(\overline{\Q}_p\) have the same \(p\)-adic
valuation \(v_p(\lambda)\), and \(v_p(\lambda)\ge0\).
\end{lemma}

\begin{proof}
The roots of \(\lambda\) are permuted transitively by
\(\Gal(\overline{\Q}_p/\Q_p)\), and the \(p\)-adic valuation is invariant
under this group, hence constant on them. Since \(M_p\) is represented by
an integral matrix on the lattice \(T_p(A)\), its characteristic roots are
algebraic integers, so their valuations are non-negative.
\end{proof}

Recall from Section~\ref{sec:setting} the subspaces $W_p$ and $U_p$ of $V_p$. Here $V_{p,\lambda}$ denotes the generalised eigenspace of $M_p$ for the irreducible factor $\lambda$. The following properties of these subspaces drive the computation. We write $E_W$ for the spectral idempotent constructed in Appendix~\ref{app:local}.

\begin{lemma}\label{sec4:act}
The following hold.
\begin{enumerate}
\item[(i)] \(W_p\) and \(U_p\) are \(M_p\)-invariant, and every
eigenvalue of \(M_p\mid W_p\) has positive \(p\)-adic valuation while every
eigenvalue of \(M_p\mid U_p\) is a \(p\)-adic unit; in particular
\(M_p\mid U_p\) is invertible and \(\det(M_p\mid U_p)\) is a \(p\)-adic
unit.
\item[(ii)] \(W_p=\{x\in V_p:\ M_p^n x\to0\ (n\to\infty)\ \text{in the }p\text{-adic
topology}\}\).
\item[(iii)] \(W_p=0\) if and only if \(p\nmid d\).
\item[(iv)] For every \(g\in G\) one has \(g(W_p)=W_p\).
\end{enumerate}
\end{lemma}

\begin{proof}
(i) By construction, every eigenvalue of $M_p\mid W_p$ has positive
valuation and every eigenvalue of $M_p\mid U_p$ is a $p$-adic unit. Hence
$M_p\mid U_p$ is invertible and has determinant in $\Z_p^\times$.

(ii) On \(W_p\) the characteristic polynomial of \(M_p\) is congruent to
\(t^{\dim W_p}\) modulo \(p\), hence some power of \(M_p\) maps the lattice
\(\Lambda_W=E_WT_p(A)\) of Appendix~\ref{app:local} into
\(p\Lambda_W\), and therefore \(M_p^n x\to0\) for every \(x\in W_p\). Conversely,
if \(x\) has a nonzero component $u\in U_p$, then all eigenvalues of
$M_p\mid U_p$ are $p$-adic units, so by (i) $\det(M_p\mid U_p)\in\Z_p^\times$. For a lattice $\Lambda_0\subset U_p$ and $r=\dim U_p$, the lattice
$\Lambda:=\sum_{i=0}^{r-1}M_p^i\Lambda_0$ satisfies $M_p\Lambda=\Lambda$, since the matrix of $M_p$ in a basis of $\Lambda$ is integral with determinant in $\Z_p^\times$. Indeed, the characteristic polynomial $\chi_U$ of $M_p\mid U_p$ lies in $\Z_p[t]$, so Cayley--Hamilton gives $M_p^r\Lambda_0\subseteq\sum_{i<r}M_p^i\Lambda_0$; the determinant is a unit, so equality follows.
With respect to $\Lambda$ the powers
$M_p^n$ and $M_p^{-n}$ are uniformly bounded. Hence $M_p^n u\to0$ with
$u\ne0$ is impossible.

(iii) The determinant of \(M_p\) on \(V_p\) is \(\pm d\), whose \(p\)-adic
valuation is \(v_p(d)\); it is the product of the eigenvalues counted with
multiplicity. Hence \(W_p=0\) if and only if all eigenvalues have
valuation \(0\), if and only if \(p\nmid d\).

(iv) Let \(g\in G\). By \(\alpha G\alpha^{-1}=G\), for every \(n\ge1\)
the automorphism \(g_n:=\alpha^n g\alpha^{-n}\) belongs to \(G\). Hence
\(M_p^n g=g_nM_p^n\) on \(V_p\). If \(x\in W_p\), then
\(M_p^n(gx)=g_n(M_p^n x)\to0\), because \(g_n\) lies in the finite group
\(G\), whose elements act as integral operators on \(T_p(A)\) (and hence are
bounded in the \(p\)-adic operator norm). Thus \(g(W_p)\subseteq W_p\);
since \(g\) is invertible and \(\dim W_p\) is finite, equality holds.
\end{proof}

\subsection{The local index lemma}
\label{sec4:lemma}

The determinant form is a statement about limits of indices of lattices. The local index lemma identifies that limit.

\begin{lemma}[Local index lemma]\label{sec4:index}
Let \(\Lambda\) be a finite free \(\Z_p\)-module, let
\(M,h\in\End(\Lambda)\) with \(M\) invertible over \(\Q_p\) and
\(hM=Mh\). Let \(W\subseteq V:=\Lambda\otimes_{\Z_p}\Q_p\) be the sum of
the generalised eigenspaces of \(M\) whose eigenvalues have positive
\(p\)-adic valuation. Then the integers
\[
I_n:=\bigl[\Lambda:h\Lambda+M^n\Lambda\bigr],\qquad n\ge0,
\]
are non-decreasing in \(n\), and
\[
\lim_{n\to\infty}I_n
 =\begin{cases}
 \bigl|\det\nolimits_{\Q_p}(h\mid W)\bigr|_p^{-1},
   & h\mid W\ \text{is invertible},\\[2pt]
 \infty, & \det\nolimits_{\Q_p}(h\mid W)=0 .
 \end{cases}
\]
The proof is given in Appendix~\ref{app:local}.
\end{lemma}

\subsection{The determinant form in the central case}
\label{sec4:central}

In the central case the local index lemma applies directly, and the determinant form follows.

\begin{lemma}[Central determinant form]\label{sec4:central-lemma}
Assume \textup{(H)} and suppose in addition that $G$ commutes with $\alpha$. Then for every \(g\in G\),
\[
\#\bigl(\Fix_A(g)\cap A[\alpha^\infty]\bigr)
=\prod_p\abs{\det\nolimits_{\Q_p}(1-g\mid W_p(\alpha))}_p^{-1}.
\]
\end{lemma}

\begin{proof}
Fix \(g\in G\). Put \(h:=1-g\in\End(A)\).
Since \(g\) commutes with \(\alpha\), the endomorphism \(h\) commutes with
\(\alpha\) and preserves every \(T_p(A)\). For each
\(n\ge1\) the endomorphism \(h\) of the finite group \(A[\alpha^n]\) has
kernel \(A[\alpha^n]\cap\Fix_A(g)\); the kernel and cokernel of an
endomorphism of a finite abelian group have the same order. Hence
\[
k_n(g)=\#\bigl(A[\alpha^n]\cap\Fix_A(g)\bigr)
=\bigl[A[\alpha^n]:(g-1)A[\alpha^n]\bigr].
\]
Further, \(A[\alpha^n]\) is the direct sum of its \(p\)-primary parts, and
for every \(p\) the $p$-primary part $A[\alpha^n]_{(p)}$ of $A[\alpha^n]$ is
isomorphic to $\alpha^{-n}T_p(A)/T_p(A)$, identified by multiplication by $\alpha^n$
with \(T_p(A)/\alpha^nT_p(A)\). Because \(g\) commutes with \(\alpha\), the action of
$g$ is compatible with this identification (indeed $g$ preserves
$\alpha^{-n}T_p(A)$). Applying Lemma~\ref{sec4:index} with
\(\Lambda=T_p(A)\), \(M=\alpha\mid T_p(A)\) and \(h=1-g\), we obtain
\[
\lim_{n\to\infty}\bigl[A[\alpha^n]_{(p)}:(g-1)A[\alpha^n]_{(p)}\bigr]
=\bigl|\det\nolimits_{\Q_p}(1-g\mid W_p)\bigr|_p^{-1},
\]
with the right-hand side read as \(\infty\) when the determinant vanishes.
Taking the product over the finitely many primes with \(W_p\ne0\) (others
contribute \(1\) by Lemma~\ref{sec4:act}(iii)) and using the
primary decomposition, we get
\[
k_\infty(g)=\#\bigl(\Fix_A(g)\cap A[\alpha^\infty]\bigr)
=\prod_p\bigl|\det\nolimits_{\Q_p}(1-g\mid W_p)\bigr|_p^{-1}.
\]
By Lemma~\ref{sec3:monotone}, $k_n(g)$ increases to $k_\infty(g)$, so this
display is the asserted formula.
\end{proof}

\subsection{Proof of the determinant form}
\label{sec4:proof}

It remains to remove the centrality assumption, which is done by replacing $\alpha$ with a centralising power.

\begin{proof}[Proof of Theorem~\ref{thm:main}(ii)]
Let \(k=\abs{\Aut(G)}\) and \(\beta=\alpha^k\) as in Lemma~\ref{sec2:beta}.
By that lemma, $\beta$ is an isogeny satisfying \eqref{eq:normalize} that
centralises $G$, and $A[\beta^\infty]=A[\alpha^\infty]$,
$W_p(\beta)=W_p(\alpha)$ for every $p$. Hence the weight of $g$ is
unchanged when $\alpha$ is replaced by $\beta$, and applying
Lemma~\ref{sec4:central-lemma} to $(A,G,\beta)$ gives
\[
\#\bigl(\Fix_A(g)\cap A[\alpha^\infty]\bigr)
=\prod_p\abs{\det\nolimits_{\Q_p}(1-g\mid W_p(\beta))}_p^{-1}
=\prod_p\abs{\det\nolimits_{\Q_p}(1-g\mid W_p(\alpha))}_p^{-1},
\]
which is the first display of Theorem~\ref{thm:main}(ii), and inserting it into
Theorem~\ref{thm:main}(i) gives the second.
\end{proof}

\subsection{Remark: quaternionic multiplication}
\label{sec4:remarks}

\begin{remark}[Quaternionic multiplication]\label{sec4:rem-qm}
The determinant form uses only $\Q_p$-linear algebra, so no further assumption on the endomorphism algebra of $A$ is needed. If \(A\) is simple with quaternionic multiplication, the finite group \(G\) acts through the embedding \(\End(A)\hookrightarrow\End(V_p)\), and \(W_p\) is \(G\)-stable by Lemma~\ref{sec4:act}(iv). Therefore the determinant form, the counting form and the limit law below apply to such abelian varieties without change. Explicit computations for them remain open (Section~\ref{sec:open}).
\end{remark}

\section{The zero case}\label{sec:zero}

Theorem~\ref{thm:main} computes $\FPP$ as the average of the weights $w_\alpha(g)$. The question here is when that average is zero. The criterion, Theorem~\ref{thm:zero}, has two equivalent forms: a fixed-point condition in $\alpha$-power torsion and a determinant condition on the active subspaces. Throughout, $V_\C := H_1(A,\C) \cong H_1(A,\Z)\otimes\C$ and by eigenvalue we mean an eigenvalue of the action of $\zeta\in G$ on $V_\C$.

\subsection{The vanishing criterion}\label{subsec:vanishing-criterion}

The determinant form of Theorem~\ref{thm:main} lets the two conditions be checked one prime at a time. Theorem~\ref{thm:zero} states them as equivalent.
\begin{theorem}[Vanishing criterion]\label{thm:zero}
Assume \textup{(H)}.  Then the following are equivalent:
\begin{enumerate}
\item[(i)] $\FPP = 0$;
\item[(ii)] $\#\bigl(\operatorname{Fix}_A(\zeta)\cap A[\alpha^\infty]\bigr) = \infty$ for every
$\zeta\in G\setminus\{1\}$;
\item[(iii)] for every $\zeta\in G\setminus\{1\}$ there exists a rational prime
$p$ with $\det\bigl(1-\zeta\mid W_p\bigr) = 0$ (that is, some $p$-adic factor of
$w_\alpha(\zeta)$ vanishes).
\end{enumerate}
For a scalar multiplier $\alpha = m\in\Z$, $|m|\ge2$, these are equivalent to
\begin{enumerate}
\item[(iv)] $D(\zeta) := \det\bigl(1-\zeta\mid H_1(A,\Z)\bigr) = 0$ for every
$\zeta\in G\setminus\{1\}$.
\end{enumerate}
\end{theorem}

\begin{proof}
By Theorem~\ref{thm:main}(i), $\FPP = |G|^{-1}\sum_{\zeta\in G} w_\alpha(\zeta)$ with
$w_\alpha(\zeta)\ge 0$, and $w_\alpha(\zeta)=0$ if and only if
$\#(\operatorname{Fix}_A(\zeta)\cap A[\alpha^\infty]) = \infty$ (with the convention
$1/\infty=0$).  This gives (i)$\Leftrightarrow$(ii).  By
Theorem~\ref{thm:main}(ii), $w_\alpha(\zeta)=\prod_p|\det(1-\zeta\mid W_p)|_p$, with a
vanishing determinant giving the factor $0$.  This vanishes if and only if
one factor vanishes, which is (iii).
For $\alpha=m$ we have $W_p=V_p$ if $p\mid m$ and $W_p=0$ otherwise, hence
$w_\alpha(\zeta)=\prod_{p\mid m}|D(\zeta)|_p$, and $|D(\zeta)|_p=0$ if and
only if $D(\zeta)=0$; this gives (iv).
\end{proof}

\begin{remark}[The positive case: comparison with Corollary~\ref{cor:pos}]
\label{rem:positive}
The zero half of the dichotomy is Theorem~\ref{thm:zero}. The positive half, $\FPP>0\iff\exists g\ne1$ with $w_\alpha(g)>0$, is Corollary~\ref{cor:pos}(a). In determinant form this reads: $\FPP>0$ if and only if there is $g\ne1$ with $\det(1-g\mid W_p)\ne0$ for every prime $p$.
\end{remark}

\begin{remark}[General isogenies: the vanishing locus depends on $\alpha$]
\label{rem:alpha-dependence}
For non-scalar $\alpha$ the condition cannot be read off $D(\zeta)$ alone.
Take $A = E_1\times E_2$ with $E_1$ and $E_2$ non-isogenous and
$G = \langle\operatorname{diag}(-1,1)\rangle\cong\mu_2$.  For $\alpha=\operatorname{diag}(2,1)$ one has
$D(-1)=0$, while $\alpha|_{E_2}=\mathrm{id}$ is an automorphism; for this datum
$\FPP=1/8$ (the limit of $1/4,5/32,17/128,\dots$).  For $\alpha=\operatorname{diag}(2,3)$
one has $\alpha|_{E_2}$ of degree $9$, and $\FPP=0$.  This is why
Theorem~\ref{thm:zero}(iii) uses the active subspace $W_p$ rather than
$H_1(A,\Z)$.
\end{remark}

\begin{proposition}[The two-dimensional case]\label{prop:dim2}
Assume \textup{(H)}.  Let $\dim A=2$, let $\alpha=m\in\Z$ with $\abs{m}\ge2$,
and let $G\le\Aut(A,0)$ be finite.  If $\FPP=0$, then $V_\C^G\ne0$.  Equivalently,
there is no I-class in dimension two (Definition~\ref{def:TI}).
\end{proposition}

\begin{proof}
Put $V_+:=\operatorname{Lie}(A)$, the tangent space of $A$ at the origin.
The Hodge decomposition is a $G$-module isomorphism
$V_\C\cong V_+\oplus\overline{V_+}$ with $\zeta$ acting as $d\zeta_0$ and
$\overline{d\zeta_0}$, so $V_\C^G=V_+^G\oplus\overline{V_+^G}$ and $V_\C^G=0$ if and
only if $V_+^G=0$.  Moreover
$\det(1-\zeta\mid V_\C)=\bigl|\det(1-d\zeta_0\mid V_+)\bigr|^2$, so $\zeta$ has
eigenvalue $1$ on $V_\C$ if and only if $d\zeta_0$ has eigenvalue $1$ on $V_+$.
The representation $\zeta\mapsto d\zeta_0$ on $V_+$ is faithful, since a
homomorphism between complex tori is determined by its differential at the origin.

Assume $\FPP=0$ and $V_\C^G=0$, so $V_+^G=0$.  By Theorem~\ref{thm:zero}(iv),
$D(\zeta)=0$ for every $\zeta\ne1$, and $D(\zeta)=\det(1-\zeta\mid V_\C)$, so
$d\zeta_0$ has eigenvalue $1$.  As $\zeta$ has finite order in characteristic
zero, $d\zeta_0$ is diagonalisable, and it is not the identity by
faithfulness.  Hence its other eigenvalue $\lambda_\zeta$ is a root of unity
different from $1$.  Let
\[
P=\frac1{\abs G}\sum_{\zeta\in G}d\zeta_0
\]
be the averaging projection on $V_+$.  Its image is $V_+^G=0$, so
$\operatorname{tr}P=0$ and
\[
0=\sum_{\zeta\in G}\operatorname{tr}_{V_+}(d\zeta_0)
 =2+\sum_{\zeta\ne1}(1+\lambda_\zeta)=\abs G+1+\sum_{\zeta\ne1}\lambda_\zeta .
\]
Thus $\bigl|\sum_{\zeta\ne1}\lambda_\zeta\bigr|=\abs G+1$; but the triangle
inequality bounds it by $\abs G-1$, a contradiction.  Hence $V_\C^G\ne0$.
\end{proof}

\subsection{Scalar multipliers: a T/I case distinction}\label{subsec:TI}

For scalar $\alpha=m\in\Z$, Theorem~\ref{thm:zero}(iv) makes the criterion explicit: $\FPP=0$ exactly when every non-identity element of $G$ has eigenvalue $1$ on $V_\C$. The distinction is then whether the common fixed space vanishes; Definition~\ref{def:TI} names the two cases.

\begin{definition}[T- and I-classes]\label{def:TI}
Let $G\le\GL(V_\C)$ be finite.  We say that $G$ is of \emph{T-class} if
$V_\C^G\ne0$, and of \emph{I-class} if $V_\C^G=0$ and every non-identity element of $G$
has eigenvalue $1$.
\end{definition}

\begin{theorem}[A case distinction for scalar multipliers]\label{thm:TI}
Let $\alpha=m\in\Z$, $|m|\ge2$, and let $G\le\Aut(A,0)$ be finite.  Then:
\begin{enumerate}
\item[(i)] $\FPP=0$ if and only if every $\zeta\in G\setminus\{1\}$ has
eigenvalue $1$ on $V_\C$.
\item[(ii)] If $G$ is of T-class, then $\FPP=0$.
\item[(iii)] An I-class example occurs for $\dim A=3$: let
$A=E_1\times E_2\times E_3$ and
let $G=V_4=\{1,a,b,c\}$ act by the sign flips
$a=(+1,-1,-1)$, $b=(-1,+1,-1)$, $c=(-1,-1,+1)$ on the three factors.  Then
$\FPP=0$ and $V_\C^G=0$.
\end{enumerate}
\end{theorem}

\begin{proof}
(i) is Theorem~\ref{thm:zero}(iv).  (ii) A common fixed vector is an
eigenvector with eigenvalue $1$ for every element, so T-class implies
$\FPP=0$ by (i).  (iii) Each non-identity element
fixes pointwise one of the three elliptic factors, hence has eigenvalue $1$
on $V_\C$; but no nonzero vector is fixed by all three elements, so $V_\C^G=0$.
\end{proof}

\begin{remark}[Status of the T/I distinction]
Every group for which the condition of Theorem~\ref{thm:TI}(i) holds is in exactly one of the two classes. The I-class is non-empty in every dimension at least three (Theorem~\ref{thm:Ifamilies}) but empty in dimension two (Proposition~\ref{prop:dim2}), so the examples start at $\dim A=3$.
\end{remark}

\subsection{I-classes: examples and infinite families}\label{subsec:Ifamilies}

The I-class is not vacuous: it occurs in every dimension at least three. The theorem below gives three families of examples and states how the proportion behaves under products.
\begin{theorem}[I-class families]\label{thm:Ifamilies}
Assume \textup{(H)}, and in (a)--(c) let $\alpha=m\in\Z$ with $|m|\ge2$.
The groups in (a)--(c) have $\FPP=0$ and $V_\C^G=0$ (hence are of I-class by Theorem~\ref{thm:TI}(i)).
Part (d) records the multiplicativity of $\FPP$ under products, where
$\FPP(A/G,\alpha)$ is the fixed-point proportion of the Latt\`es map attached to $(A,G,\alpha)$.
\begin{enumerate}
\item[(a)] (\emph{$V_4$ block-flip family}) Let $A=E_1\times\cdots\times E_g$ with
$\dim A=g\ge3$ and $E_1,\dots,E_g$ elliptic curves.  Partition
$\{1,\dots,g\}$ into three non-empty classes $C_0,C_1,C_2$ and define
$a,b,c\in\Aut(E_1\times\cdots\times E_g)$ by
\[
a = \mathrm{diag}(+1 \text{ on } C_0,\,-1 \text{ on } C_1\cup C_2),\quad
b = \mathrm{diag}(-1 \text{ on } C_0\cup C_2,\,+1 \text{ on } C_1),\quad
c=ab .
\]
Then $G=\{1,a,b,c\}\cong V_4$ satisfies $\FPP=0$ and $V_\C^G=0$.
\item[(b)] (\emph{Even-weight sign groups}) Let $E$ be an elliptic curve, let
$g=\dim A\ge3$ be odd, put $A=E^g$, and let
$G=\{\varepsilon\in\{\pm1\}^g : \prod_{i=1}^g\varepsilon_i = 1\}
\cong(\Z/2)^{g-1}$ act diagonally on $E^g$.  Then $\FPP=0$ and $V_\C^G=0$.
\item[(c)] (\emph{Tetrahedral and octahedral groups}) Let $E$ be a non-CM
elliptic curve, $A=E^3$, and identify $\Aut(E^3,0)=\GL_3(\Z)$.  Let $G\in\{A_4,S_4\}$ be the
tetrahedral or octahedral rotation group, respectively, embedded in
$\mathrm{SO}(3)\cap\GL_3(\Z)$.  Then $\FPP=0$ and $V_\C^G=0$.
\item[(d)] (\emph{Products}) For the product, let $A=A_1\times A_2$ and
$\alpha=\alpha_1\times\alpha_2$ with $\alpha_i$ an isogeny of degree at least
$2$ normalising $G_i\le\Aut(A_i,0)$, and let
$G=G_1\times G_2$ act factorwise.  Then, for $(u_1,u_2)\in G=G_1\times G_2$,
$w_\alpha(u_1,u_2)=w_{\alpha_1}(u_1)w_{\alpha_2}(u_2)$ and
$\FPP(A/G,\alpha)=\FPP(A_1/G_1,\alpha_1)\FPP(A_2/G_2,\alpha_2)$.  In particular, if one factor has
zero proportion, so does the product.
\end{enumerate}
\end{theorem}

\begin{proof}
For (a)--(c) it suffices to find an eigenvalue $1$ on $V_\C$ for every
non-identity element, since then $D(\zeta)$ vanishes and $\FPP=0$ follows from
Theorem~\ref{thm:zero}(iv).

(a) Here $V_\C=\bigoplus_{j=1}^gH_1(E_j,\C)$ with $\dim_\C H_1(E_j,\C)=2$, and
the signs of $a,b,c$ on $(C_0,C_1,C_2)$ are $(+1,-1,-1)$, $(-1,+1,-1)$ and
$(-1,-1,+1)$.  A non-identity element acts by $+1$ on one class $C_i$ and by
$-1$ on the other two, so it has eigenvalue $1$ with multiplicity
$2\abs{C_i}\ge2$.  A vector $v=\sum_jv_j$ fixed by $a$ and $b$ satisfies
$v_j=0$ on $C_1\cup C_2$ (from $a$) and on $C_0\cup C_2$ (from $b$), so
$v=0$.

(b) Here $A=E^g$, $V_\C=\bigoplus_{i=1}^gH_1(E_i,\C)$ with $E_1=\cdots=E_g=E$,
and $\varepsilon\in G$ acts on the $i$-th summand as $\varepsilon_i$.  A
non-identity element is not the vector with all entries $-1$, because
$(-1)^g=-1$ for odd $g$.  Hence some $\varepsilon_i=+1$, and $\varepsilon$
has eigenvalue $1$ on $H_1(E_i,\C)$.  If $v=(v_1,\dots,v_g)$ is fixed by $G$,
fix $i$, choose $j\ne i$ (as $g\ge3$), and take $\varepsilon\in G$ equal to
$-1$ at $i$ and $j$ and $+1$ elsewhere.  Then $(\varepsilon v)_i=-v_i$, so
$v_i=0$ and $v=0$.

(c) Since $E$ is non-CM, $\End(E)=\Z$, hence $\End(E^3)=M_3(\Z)$ and
$\Aut(E^3,0)=\GL_3(\Z)$.  Then $V_\C=H_1(E^3,\C)\cong H_1(E,\C)\otimes_\Z\Z^3$,
with $G$ acting trivially on $H_1(E,\C)$ and by its defining integral matrix
on $\Z^3$, so each eigenvalue of $\zeta\in G$ on $\C^3$ occurs on $V_\C$ with
multiplicity two.  A non-identity element of $G\le\mathrm{SO}(3)$ is a
nontrivial rotation, hence has eigenvalue $1$ on $\C^3$ and on $V_\C$.  Also
$V_\C^G\cong H_1(E,\C)\otimes(\C^3)^G$, so it suffices that $(\C^3)^G=0$.

The rotation representation of $A_4$ is irreducible: $A_4$ has irreducible
degrees $1,1,1,3$, and this one is
faithful with non-abelian image, so it is not a sum of one-dimensional
representations.  An irreducible representation of dimension greater than
one has no nonzero invariant vector, so $(\C^3)^{A_4}=0$, and then
$(\C^3)^G=0$ for $G=A_4$ and for $G=S_4$ as $A_4\le S_4$.  For $S_4$ the rotation
representation is $\mathrm{std}\otimes\mathrm{sign}$ rather than
$\mathrm{std}$.  It is irreducible, as its restriction to $A_4$ is, and the
only irreducible representations of $S_4$ of degree three are
$\mathrm{std}$ and $\mathrm{std}\otimes\mathrm{sign}$.  It is not
$\mathrm{std}$, because a $4$-cycle has eigenvalues $i,-1,-i$ on
$\mathrm{std}$, whereas a $90^\circ$ rotation has eigenvalue $1$.

(d) For $u=(u_1,u_2)\in G_1\times G_2$ we have
$\Fix_A(u)=\Fix_{A_1}(u_1)\times\Fix_{A_2}(u_2)$ and
$A[\alpha^\infty]=A_1[\alpha_1^\infty]\times A_2[\alpha_2^\infty]$, so the
intersection is a product and the cardinalities multiply (an infinite factor
giving $\infty$, with $1/\infty=0$).  Hence
$w_\alpha(u_1,u_2)=w_{\alpha_1}(u_1)w_{\alpha_2}(u_2)$.  Averaging over $G$
and applying Theorem~\ref{thm:main}(i),
\[
\FPP(A/G,\alpha)=\frac1{\abs{G_1}\abs{G_2}}\sum_{u_1\in G_1}\sum_{u_2\in G_2}
w_{\alpha_1}(u_1)w_{\alpha_2}(u_2)=\FPP(A_1/G_1,\alpha_1)\FPP(A_2/G_2,\alpha_2),
\]
so if one factor has zero proportion, so does the product.
\end{proof}

\section{The fixed-point process}\label{sec:process}

For $n\ge1$ let $Y_n\colon G_\infty\to\Z_{\ge0}$ be defined by $Y_n(g)=\#\operatorname{Fix}(\pi_n(g)\mid L_n)$, the number of vertices fixed at level $n$. Then $\FPP=\Prob[Y_n>0\ \text{eventually}]$. The finite-level structure of Section~\ref{sec:finite-level} gives $Y_n$ exactly, and Theorem~\ref{thm:main} evaluates the resulting average. The two results below identify its behaviour along the tower; Remark~\ref{rem:LL-scope} records the scope of the limit law.

\subsection{The martingale property}\label{subsec:martingale}

The translations provide the transitivity on children that the martingale criterion asks for. The next theorem states the consequence, and Section~\ref{subsec:limit-law} passes to the limit.
\begin{theorem}[Martingale property]\label{thm:martingale}
Assume only \textup{(H)}.  Then
$G_\infty(K,f,t)$ is a martingale group, i.e. its fixed-point process is a
martingale \cite[Theorem~5.7]{BJKL}:
for every $n\ge1$ and every vertex $v$ at level $n-1$, the kernel
$N_n=\ker(G_n\to G_{n-1})$ acts transitively on the children $v^*$ of $v$.
Consequently $\{Y_n\}$ is a martingale with respect to the natural
filtration; $Y_n$ converges almost surely to an integer-valued random
variable $Y_\infty$, and $\FPP=\Prob[Y_\infty>0]$.
\end{theorem}

\begin{proof}
By Lemma~\ref{sec3:affine}, at level $n$ an element $\sigma$ with
$\sigma\mid K(A)=\zeta$ acts by $h\mapsto\gamma_n^{-1}(h+c)$, where
$\gamma_n=\alpha^{-n}\zeta\alpha^n$.  For $\zeta=1$ this is the translation
$h\mapsto h+c$.  A translation acts trivially on level $n$ only for $c=0$;
hence the translation subgroup maps injectively into $G_n$ and $G_n$ contains
the images of all translations by $A[\alpha^n]$.  In particular $G_n$ acts
transitively on the vertices of level $n$, so $G_\infty$ is level-transitive.

Let $c\in A[\alpha]=\ker\alpha$.  The translation by $c$ sends a level-$n$
vertex $h$ to $h+c$, whose parent is $\alpha(h+c)=\alpha h$; so it fixes every
level-$(n-1)$ vertex and its image lies in $N_n$.  The children of $v$ are
the points $h\in A[\alpha^n]$ with $\alpha h=v$, and they form a torsor under
$A[\alpha]$: if $h_1,h_2\in v^*$, then $\alpha(h_1-h_2)=0$, so
$h_1-h_2\in A[\alpha]$, and the translation by $h_2-h_1$ maps $h_1$ to
$h_2$.  Thus $N_n$ acts transitively on $v^*$.  The criterion
\cite[Theorem 5.7]{BJKL} (equivalently \cite[Theorem 3.5]{HZ}) now gives
that $\{Y_n\}$ is a martingale; the martingale convergence theorem and the
identity $\FPP=\Prob\big[Y_n>0\ \text{for all }n\big]$ finish the proof.
\end{proof}

\subsection{The limit law}\label{subsec:limit-law}

The martingale gives convergence. The exact law of $Y_\infty$ is read off from the finite-level Haar measure: the linear parts are the conjugates $\gamma_n=\alpha^{-n}\zeta\alpha^n$ of $\zeta$, and the determinant form enters through $w_\alpha(\zeta)=\kappa(\zeta)^{-1}$ (Section~\ref{sec:setting}).

\begin{lemma}[Independence of the level]\label{sec6:residue}
Assume \textup{(H)}.  Let $\zeta\in G$ and
$\gamma_r=\alpha^{-r}\zeta\alpha^r$ for $r\ge0$.  Then for every prime $p$
and every $r\ge0$,
$\det\bigl(1-\gamma_{r+1}\mid W_p(\alpha)\bigr)=\det\bigl(1-\gamma_r\mid W_p(\alpha)\bigr)$.
\end{lemma}

\begin{proof}
The determinants are defined since $W_p(\alpha)$ is $G$-stable
(Lemma~\ref{sec4:act}(iv)).  From $\gamma_{r+1}=\alpha^{-1}\gamma_r\alpha$
we get $\alpha\gamma_{r+1}=\gamma_r\alpha$ in $\End(A)$, hence
$\rho(\gamma_{r+1})=M_p^{-1}\rho(\gamma_r)M_p$ on $V_p$, where
$M_p=\alpha\mid V_p$ and $\rho$ is the action of $G$.  Every eigenvalue of
$M_p\mid W_p(\alpha)$ is nonzero, so $M_p^{\pm1}$ preserves
$W_p(\alpha)$.  Hence
$1-\gamma_{r+1}=M_p^{-1}(1-\gamma_r)M_p$ on $W_p(\alpha)$, and
multiplicativity of the determinant gives
$\det(1-\gamma_{r+1}\mid W_p(\alpha))=\det(1-\gamma_r\mid W_p(\alpha))$.
\qedhere
\end{proof}

\begin{theorem}[Limit law]\label{thm:limitlaw}
Assume \textup{(H)}, and let $\kappa(\zeta)$ be as in Section~\ref{sec:setting}.
Let $\zeta$ be uniform on $G$ and $c$ be Haar-random on
$T_\alpha(A)=\varprojlim_n A[\alpha^n]$, independent of $\zeta$.  Put
$\gamma_n=\alpha^{-n}\zeta\alpha^n$,
$k_n(g)=\#\bigl(A[\alpha^n]\cap\Fix_A(g)\bigr)$, and $\tilde Y_n=k_n(\gamma_n)$
if $c_n\in(\gamma_n-1)A[\alpha^n]$, $\tilde Y_n=0$ otherwise, where $c_n$ is the image of $c$ in $A[\alpha^n]$.  Then
$(\tilde Y_n)$ has the same law as $(Y_n)$ for a Haar-random element of
$G_\infty$.  Almost surely it equals $\kappa(\zeta)$ if
$c_n\in(\gamma_n-1)A[\alpha^n]$ for every $n$, and $0$ otherwise.  On the
null set where $\kappa(\zeta)=\infty$ and $c_n\in(\gamma_n-1)A[\alpha^n]$ for every $n$, the limit is $+\infty$.
With $1/\infty=0$, conditionally on $\zeta$,
\[
\Prob\left[\lim_n\tilde Y_n=\kappa(\zeta)\;\middle|\;\zeta\right]
=\frac{1}{\kappa(\zeta)},
\qquad
\Prob\left[\lim_n\tilde Y_n=\kappa(\zeta)\right]
=\frac1{\abs{G}}\sum_{\zeta\in G}\frac1{\kappa(\zeta)}
=\FPP.
\]
For scalar $\alpha=m\in\Z$, $\abs m\ge2$, and $\zeta\ne1$ with
$D(\zeta)\ne0$, $\kappa(\zeta)=\rho_m(\zeta)^{-1}$ with
$\rho_m(\zeta)=\prod_{p\mid m}\abs{D(\zeta)}_p$.  Hence for each positive
value $r$,
\[
\Prob\left[\lim_n\tilde Y_n=r\right]
=\frac1{\abs G}\sum_{\substack{\zeta\in G\\ \rho_m(\zeta)^{-1}=r}}\rho_m(\zeta),
\]
the sum being over $1\ne\zeta$ with $D(\zeta)\ne0$ (grouping equal values).
\end{theorem}

\begin{proof}
\textbf{Step 1.} By Theorem~\ref{sec3:galois}, $\Gamma_n\cong A[\alpha^n]\rtimes G$ has order $d^n\abs G$.  By Lemma~\ref{sec3:affine}, $\sigma\mapsto(\zeta,c)$ is a bijection $\Gamma_n\to G\times A[\alpha^n]$, so Haar measure on $\Gamma_n$ is the uniform product on $G\times A[\alpha^n]$.  By Theorem~\ref{sec3:galois} and the first paragraph of the proof of Proposition~\ref{sec3:fppn-formula}, $G_n=\Gamma_n/H_n$ with $H_n$ acting trivially on the level-$n$ vertices, so $H_n$-invariant proportions agree in $\Gamma_n$ and $G_n$.  Since $G_\infty\to G_n$ is surjective, the law of $Y_n$ is that of the functional of independent uniform $\zeta\in G$ and $c_n\in A[\alpha^n]$.

\textbf{Step 2.} The inclusions $K(y_n)\subseteq K(y_{n+1})$ give surjective $\Gamma_{n+1}\to\Gamma_n$, so Haar measure on $\Gamma_\infty=\varprojlim_n\Gamma_n$ pushes forward to Haar measure on every $\Gamma_n$.  For a cylinder set defined by conditions at levels $n_1<\cdots<n_N$, let $C\subseteq\Gamma_{n_N}$ be the set of elements realising them.  Then $C$ is a union of cosets of $H_{n_N}$, so its proportion is $\abs C/\abs{\Gamma_{n_N}}$, and this equals the probability of the cylinder under Haar measure on $G_\infty$, on $G_{n_N}$ and on $\Gamma_{n_N}$.  By Step 1 it also equals the probability of the model pair $(\zeta,c_{n_N})$ with the lower-level parameters $\alpha^{n_N-n_i}c_{n_N}$.  Hence the finite-dimensional distributions agree.  In particular, no identification of $G_\infty$ with $\Gamma_\infty$, and no description of $G_\infty$ as a subgroup of $T_\alpha(A)\rtimes G$, is required.

\textbf{Step 3.} Put $S_n:=(\gamma_n-1)A[\alpha^n]$ and $E_n:=\{c:c_n\in S_n\}$.  From $\gamma_{n+1}=\alpha^{-1}\gamma_n\alpha$ and $\alpha A[\alpha^{n+1}]=A[\alpha^n]$,
\[
\alpha S_{n+1}=\alpha(\gamma_{n+1}-1)A[\alpha^{n+1}]=(\gamma_n-1)A[\alpha^n]=S_n,
\]
so $E_{n+1}\subseteq E_n$ since $c_n=\alpha c_{n+1}$.  As $c\mapsto c_n$ pushes Haar measure to the uniform measure and $\gamma_n-1$ has kernel $A[\alpha^n]\cap\Fix_A(\gamma_n)$ and image $S_n$, $\mu_n:=\Prob\big[E_n\mid\zeta\big]=1/k_n(\gamma_n)$ is non-increasing.

\textbf{Step 4.} By Lemma~\ref{sec3:affine}, $\tilde Y_n=k_n(\gamma_n)$ on $E_n$ and $0$ outside.  If $c\notin\bigcap_nE_n$, then $c\notin E_{n_0}$ for some $n_0$, so $\tilde Y_n=0$ for $n\ge n_0$ and $\tilde Y_n\to0$.

\textbf{Step 5.} By Lemma~\ref{lem:E1}(iii), the conjugation $c_\alpha\colon G\to G$ is an automorphism of $G$.  Put $k=\abs{\Aut(G)}$, so $c_\alpha^k=\mathrm{id}$ and $\gamma_{n+k}=\gamma_n$ for all $n$.  Step 3 makes $k_n(\gamma_n)$ non-decreasing.  For $r\in\{0,\dots,k-1\}$ and $n\equiv r\pmod k$ we have $\gamma_n=\gamma_r$, so Lemma~\ref{sec3:monotone}, Theorem~\ref{thm:main}(ii) and Lemma~\ref{sec6:residue} give
\[
\lim_{\substack{n\to\infty\\ n\equiv r}}k_n(\gamma_n)=k_\infty(\gamma_r)=\prod_p\abs{\det(1-\gamma_r\mid W_p(\alpha))}_p^{-1}=\prod_p\abs{\det(1-\zeta\mid W_p(\alpha))}_p^{-1}=\kappa(\zeta).
\]
All $k$ residue classes give the same limit, so $k_n(\gamma_n)\to \kappa(\zeta)$ monotonically and $\mu_n\to1/\kappa(\zeta)$ with $1/\infty=0$.  By continuity from above, $\Prob\big[\bigcap_nE_n\mid\zeta\big]=1/\kappa(\zeta)$, and on $\bigcap_nE_n$ we have $\tilde Y_n=k_n(\gamma_n)\to \kappa(\zeta)$.  If $\kappa(\zeta)=\infty$, the set $\bigcap_nE_n$ is null and the limit is $+\infty$ on it.

\textbf{Step 6.} Since $\zeta$ is uniform on $G$ and independent of $c$,
\[
\Prob\left[\lim_n\tilde Y_n=\kappa(\zeta)\right]=\frac1{\abs G}\sum_{\zeta\in G}\Prob\left[\bigcap_nE_n\;\middle|\;\zeta\right]=\frac1{\abs G}\sum_{\zeta\in G}\frac1{\kappa(\zeta)}.
\]
Since $\kappa(\zeta)\ge1$, the limit is positive exactly on $\bigcap_nE_n$.  By Step 2 both processes have the same law, so $\FPP=\Prob\bigl[\bigcap_nE_n\bigr]$ and
\[
\FPP=\frac1{\abs G}\sum_{\zeta\in G}\frac1{\kappa(\zeta)},
\]
which is the counting form of Theorem~\ref{thm:main}(i) as $w_\alpha(\zeta)=\kappa(\zeta)^{-1}$.  For scalar $\alpha=m$, $W_p(\alpha)=V_p$ for $p\mid m$ and $0$ otherwise (Lemma~\ref{sec4:act}(iii)), so $\kappa(\zeta)=\rho_m(\zeta)^{-1}$ if $D(\zeta)\ne0$.
\end{proof}

\begin{example}[An exact limit law]\label{ex:limit-laws}
For $A=E^2$ with $E:y^2=x^3-x$, $G=\langle-I\rangle\cong\mu_2$ and $m=2$,
one has $Y_\infty=16$ with probability $1/32$ and $\FPP=1/32$.  All values
follow from the determinant computation in Theorem~\ref{thm:main}(ii).
\end{example}

\begin{remark}[Scope of the limit law]\label{rem:LL-scope}
Theorem~\ref{thm:limitlaw} holds under \textup{(H)} alone. The proof uses the finite-level Haar measure, the nesting of the events $E_n$ and the invariance from Lemma~\ref{sec6:residue}, together with the results of Sections~\ref{sec:setting}--\ref{sec:det}. The linear part at level $n$ is the conjugate $\gamma_n$ rather than $\zeta$.
\end{remark}

\subsection{Geometric instances of the value range for Question 8}
\label{subsec:Q8}

Question~8 of \cite{FARQ8} asks which values the limit can take. For abstract level-transitive groups the possible values are exactly $\N\cup\{0\}$ (\cite[{\S}2]{FARQ8}), and \cite[Theorem~9]{FARQ8} realises the value $d$. The corollary below shows that the value $d$, and more generally the values $\rho_m(\zeta)^{-1}$, occur for geometric iterated Galois groups of Latt\`es maps.

\begin{corollary}[Geometric instance of the value range]\label{cor:Q8}
Assume \textup{(H)} and let $\alpha=m\in\Z$, $|m|\ge2$.  Then
$G_\infty(K,f,t)$ is level-transitive and its fixed-point process converges
almost surely; the positive values of $Y_\infty$ occurring with positive
probability are exactly $\rho_m(\zeta)^{-1}$ for $1\ne\zeta\in G$ with
$D(\zeta)\ne0$, and their probabilities are $\rho_m(\zeta)/|G|$
(grouping equal values).  In particular,
for $A=E^2$ ($E:y^2=x^3-x$), $G=\mu_4^2$ and $m=2$ (so $d=m^{2\dim A}=16$), the
positive limit values are
\[
r\in\{4,\,8,\,16\},\qquad
\Prob[r=4]=1/16,\quad \Prob[r=8]=1/32,\quad \Prob[r=16]=1/256,
\]
and $\FPP=25/256$.
\end{corollary}

\begin{proof}
Level-transitivity follows from the finite-level structure: $G_n$ contains
the translations $A[m^n]$.  The limit law is Theorem~\ref{thm:limitlaw};
the values follow from the multiset of $\det(1-\zeta\mid W_2)$ over $G=\mu_4^2$, namely
$\{0$ (seven times), $4$ (four times), $8$ (four times), $16$ (once)$\}$.
\end{proof}

\begin{remark}[The one-dimensional shadow]\label{rem:g1shadow}
Limit values outside $\{0,1,2,d\}$ already occur for $\dim A=1$: for an elliptic
curve $E$ with $G=\mu_3\le\Aut(E)$, $\alpha=m=3$ (so $d=9$), and a primitive
cube root $\omega$ of unity, $\#\ker(1-\zeta\mid E[3^n])=3$ for
$\zeta=\omega,\omega^2$ and every $n$, so Theorem~\ref{thm:limitlaw} gives
$Y_\infty=3$ with probability $2/9$.  Thus $3\notin\{0,1,2,9\}$ occurs for
$\dim A=1$.

For $E$ with $G=\mu_6\le\Aut(E)$ and $m=6$ (so $d=36$), the nonzero
determinants of the non-identity elements are $1,1,3,3,4$, so
Theorem~\ref{thm:limitlaw} gives positive limit values $1,3,4$ with
probabilities $1/3,1/9,1/24$, all of them below $d$.  All values are computed
from Corollary~\ref{cor:scalar} together with Theorem~\ref{thm:limitlaw}.
\end{remark}

\section{The one-dimensional case: ten exact values for scalar multipliers}\label{sec:onedim}

In dimension one the general theorems of Sections~\ref{sec:setting}--\ref{sec:det} become explicit. The formula yields ten values in closed form for scalar multipliers, and the Bridy--Jones--Kelsey--Lodge bound is attained (Corollary~\ref{cor:BJKL}). The exceptional case of the Fari\~na-Asategui--Radi fixed-point proportion conjecture is settled over separably closed fields of characteristic $0$, and the eight quadratic Latt\`es classes are determined.

\subsection{The ten values for scalar multipliers}

Let $A=E$ be an elliptic curve over $K$, so that $\dim A=1$, let $\alpha=m\in\Z$ with $\abs{m}\ge2$, and let $G\le\Aut(E,0)$ be nontrivial. The group $\Aut(E,0)$ is cyclic of order $2$, $4$ or $6$, so $G$ is one of $\mu_2,\mu_3,\mu_4,\mu_6$. If $\abs{G}\ge3$ then $E$ has complex multiplication: the inclusion $\mu_3\le\Aut(E,0)$ forces $j(E)=0$, and $\mu_4\le\Aut(E,0)$ forces $j(E)=1728$. Only $\mu_2$ occurs on every elliptic curve, and the values for $G=\mu_2$ below are the same with and without complex multiplication. For a non-scalar isogeny $\alpha$ the value is given by Theorem~\ref{thm:main}(ii).

Fix $\zeta\in G$ with $\zeta\ne1$ and put $\abs{x}_m:=\prod_{p\mid m}\abs{x}_p$ for $x\in\Q$. By Theorem~\ref{thm:main}(ii), the number
\[
k_\infty(\zeta):=\#\bigl(\Fix_E(\zeta)\cap E[m^\infty]\bigr)
=\prod_{p\mid m}\abs{\det(1-\zeta\mid H_1(E,\Z))}_p^{-1}
=\abs{\det(1-\zeta\mid H_1(E,\Z))}_m^{-1}
\]
is finite. It is also the index of $(\zeta-1)E[m^\infty]$ in $E[m^\infty]$: for every $n$ the kernel and the cokernel of $1-\zeta$ on the finite group $E[m^n]$ have the same order, and both orders stabilise as $n$ grows. The determinant is computed as follows. If $E$ has complex multiplication, then $\End(E)\cong\OO$ for an order $\OO$ in an imaginary quadratic field $F$, and $H_1(E,\Z)\otimes_\Z\Q\cong F$ as $F$-vector spaces, so that $\det(1-\zeta\mid H_1(E,\Z))=N_{F/\Q}(1-\zeta)$. If $E$ has no complex multiplication, then $\End(E)=\Z$, $G=\mu_2$ and $\zeta=-1$, so that $\det(1-\zeta\mid H_1(E,\Z))=\det(2I_2)=4$. The resulting determinant values are used in the proof of Corollary~\ref{cor:values1}.

\begin{corollary}\label{cor:values1}
Assume \textup{(H)} with $A=E$ an elliptic curve, $\alpha=m\in\Z$, $\abs{m}\ge2$, and $G\le\Aut(E,0)$ nontrivial.
\begin{enumerate}
\item[(i)] If $G=\mu_2$, then $\FPP=1/8$ when $2\mid m$ and $\FPP=1/2$ when $2\nmid m$.
\item[(ii)] If $G=\mu_4$, then $\FPP=5/16$ when $2\mid m$ and $\FPP=3/4$ when $2\nmid m$.
\item[(iii)] If $G=\mu_3$, then $\FPP=2/9$ when $3\mid m$ and $\FPP=2/3$ when $3\nmid m$.
\item[(iv)] If $G=\mu_6$, then
\[
\FPP=\begin{cases}
35/72, & 6\mid m,\\
17/24, & 2\mid m,\ 3\nmid m,\\
11/18, & 2\nmid m,\ 3\mid m,\\
5/6, & (m,6)=1
\end{cases}.
\]
\end{enumerate}
In particular, for scalar multipliers, the proportion is positive for every one-dimensional Latt\`es map with nontrivial $G\le\Aut(E,0)$, but it is not determined by the group $G$ alone: it also depends on which primes divide $m$.
\end{corollary}

\begin{proof}
By Corollary~\ref{cor:scalar},
\[
\FPP=\frac1{\abs{G}}\sum_{1\ne\zeta\in G}\abs{\det(1-\zeta\mid H_1(E,\Z))}_m .
\]
The determinant values $\det(1-\zeta\mid H_1(E,\Z))$ are
\[
\begin{cases}
4, & G=\mu_2,\ \zeta=-1,\\
2,\ 4,\ 2, & G=\mu_4,\ \zeta=i,-1,-i,\\
3, & G=\mu_3,\ \zeta=\omega,\omega^2,\\
1,\ 1,\ 3,\ 3,\ 4, & G=\mu_6,\ \zeta=\zeta_6,\zeta_6^{-1},\omega,\omega^2,-1
\end{cases},
\]
where $\omega$ is a primitive cube root of unity and $\zeta_6=e^{\pi i/3}$. The index of $(\zeta-1)E[m^\infty]$ in $E[m^\infty]$ equals $1/\abs{\det(1-\zeta\mid H_1(E,\Z))}_m$: it is $4$ or $1$ for determinant $4$ according as $2\mid m$ or not, $2$ or $1$ for determinant $2$ according as $2\mid m$ or not, $3$ or $1$ for determinant $3$ according as $3\mid m$ or not, and $1$ for determinant $1$.
For $G=\mu_2$ the only non-identity element is $-1$, with determinant $4$; hence $\abs{4}_m=1/4$ when $2\mid m$ and $\abs{4}_m=1$ otherwise, and the average is $(1/2)(1/4)=1/8$ in the first case and $(1/2)(1)=1/2$ in the second.

For $G=\mu_4$ the non-identity elements are $i,-1,-i$, with determinants $2,4,2$. Hence these determinants have $\abs{\cdot}_m$ equal to $1/2,1/4,1/2$ when $2\mid m$ and to $1,1,1$ otherwise, so the average is $(1/4)(1/2+1/4+1/2)=5/16$ in the first case and $(1/4)(1+1+1)=3/4$ in the second.

For $G=\mu_3$ the non-identity elements are $\omega,\omega^2$, with determinant $3$. Hence $\abs{3}_m=1/3$ when $3\mid m$ and $\abs{3}_m=1$ otherwise, so the average is $(1/3)(1/3+1/3)=2/9$ in the first case and $(1/3)(1+1)=2/3$ in the second.

For $G=\mu_6$ the non-identity elements are $\zeta_6,\zeta_6^{-1}$ with determinant $1$, $\omega,\omega^2$ with determinant $3$, and $-1$ with determinant $4$. Hence
\[
\FPP=\frac16\bigl(2+2\abs{3}_m+\abs{4}_m\bigr).
\]
If $6\mid m$, this is $\frac16(2+\frac23+\frac14)=\frac{35}{72}$. If $2\mid m$ and $3\nmid m$, it is $\frac16(2+2+\frac14)=\frac{17}{24}$. If $2\nmid m$ and $3\mid m$, it is $\frac16(2+\frac23+1)=\frac{11}{18}$. If $(m,6)=1$, it is $\frac16(2+2+1)=\frac56$.
\end{proof}

In dimension one, the upper bound $(\abs G-1)/\abs G$ of Corollary~\ref{cor:pos}(b) is attained exactly for $\mu_2$ and $\mu_4$ when $2\nmid m$, for $\mu_3$ when $3\nmid m$, and for $\mu_6$ when $(m,6)=1$.

The values are collected in Table~\ref{tab:tenvalues}.

\begin{table}[tbp]
\centering
\small
\begin{tabular}{@{}lcccc@{}}
\toprule
$G$ & $2\mid m,\ 3\nmid m$ & $2\nmid m,\ 3\mid m$ & $6\mid m$ & $(m,6)=1$ \\
\midrule
$\mu_2$ & $1/8$ & $1/2$ & $1/8$ & $1/2$ \\
$\mu_3$ & $2/3$ & $2/9$ & $2/9$ & $2/3$ \\
$\mu_4$ & $5/16$ & $3/4$ & $5/16$ & $3/4$ \\
$\mu_6$ & $17/24$ & $11/18$ & $35/72$ & $5/6$ \\
\bottomrule
\end{tabular}
\caption{The ten values of Corollary~\ref{cor:values1} for scalar multipliers. The four columns are the four divisibility cases. The ten distinct values appear with repetitions among the sixteen entries. For $\mu_2$ and $\mu_4$ the value depends only on $2\mid m$, and for $\mu_3$ only on $3\mid m$, so those entries repeat across the columns.}
\label{tab:tenvalues}
\end{table}

\subsection{Attainment of the bound \texorpdfstring{$1/8$}{1/8}}

The bound is due to Bridy--Jones--Kelsey--Lodge: \cite[Theorem~2.5]{BJKL} gives the lower bound $1/8$ for the periodic-point proportion of a class of quadratic Latt\`es maps. The corollary below evaluates the geometric proportion for one map of that class.
\begin{corollary}[The bound $1/8$ is attained]\label{cor:BJKL}
Let $E:y^2=x^3+x$ over $K$, let $G=\{\pm1\}$ and let $\alpha=1+i$, so that $\End(E)=\Z[i]$. Then $\deg f=N_{\Q(i)/\Q}(1+i)=2$, the map $f$ is linearly conjugate to the Latt\`es map $(x^2+a)/(x^2-a)$ with $a^2=-1$, and
\[
\FPP(G_\infty)=\frac18 .
\]
Consequently the bound $1/8$ in \cite[Theorem 2.5]{BJKL} is attained, namely by the reductions of $f$ modulo the split primes $p\equiv1\pmod4$ of good reduction.
\end{corollary}

\begin{proof}
The value $1/8$ is computed below from Theorem~\ref{thm:main}(ii) without any input from \cite{BJKL}. The comparison with the limit inferior that follows uses \cite[Corollary 3.5]{BJKL}, \cite[Theorem 3.9, (3.10) and \S4]{BJKL} and \cite[Theorem 5.19]{R}; the sharpness statement for that limit inferior therefore combines the computation with these cited results.
The eigenvalues of $\alpha$ on $V_2=T_2(E)\otimes_{\Z_2}\Q_2$ are $1+i$ and $1-i$, and both have $2$-adic valuation $1/2$, because $(1+i)^2=2i$ has valuation $1$. Hence $W_2=V_2$. For $\zeta=-1$ we have $\det(1-\zeta\mid W_2)=\det(2\mid V_2)=4$, and $\abs{4}_2=1/4$. Since $\deg\alpha=2$, Lemma~\ref{sec4:act}(iii) gives $W_p=0$ for every prime $p\ne2$. Theorem~\ref{thm:main}(ii) therefore gives
\[
\FPP(G_\infty)=\frac12\cdot\abs{4}_2=\frac12\cdot\frac14=\frac18 .
\]
For the comparison with \cite{BJKL}, the data $(E,G,\alpha)$ are defined over $\Q(i)$; let $f_0$ be the corresponding map over $\Q(i)$, and let $f_\C$ be its base change to $\C$. By \cite[Theorem 5.19]{R}, whose hypotheses hold over $\C$, the geometric iterated Galois group of $f_\C$ is its iterated monodromy group $\operatorname{IMG}(f_\C)$; the same computation over $\C$ gives
\[
 \FPP(\operatorname{IMG}(f_\C))=\FPP(G_\infty(\C,f_\C,t))=1/8 .
\]
In the $x$-coordinate the map $[\alpha]$ descends to $f_0(x)=-i(x^2+1)/(2x)$. Put $g(z)=(z^2+i)/(z^2-i)$ and $M(z)=-i(z-1)/(z+1)$. Then $M(f_0(x))=g(M(x))$ for every $x$, so $f_0$ is linearly conjugate to $g$, and hence to $(x^2+a)/(x^2-a)$ with $a^2=-1$.

Choose a split rational prime $p\equiv1\pmod4$ of good reduction; such primes exist by Dirichlet's theorem. Let $f_p$ be the reduction of $f_0$ modulo a prime above $p$. Then $f_p$ is defined over $\Fp_p$, and by \cite[Corollary 3.5]{BJKL}, for this fixed split $p$,
\[
 \liminf_{n\to\infty}\frac{\#\operatorname{Per}(f_p,\PP^1(\Fp_{p^n}))}{p^n+1}\le\FPP\bigl(\operatorname{pgIMG}(f_p)/\Fp_p\bigr),
\]
where $\operatorname{Per}$ denotes the set of periodic points and $\operatorname{pgIMG}$ the pro-geometric IMG (see \cite[\S3]{BJKL}).
For $p$ odd, the critical points $\pm1$ and their postcritical orbits $\pm i\to0\to\infty$ reduce injectively modulo any prime above $p$, so $f_p$ and $f_\C$ have the same ramification portrait; by Grothendieck's specialisation theorem for tame fundamental groups (the argument of \cite[Theorem~3.9]{BJKL}) the tree actions are conjugate, hence \cite[(3.10)]{BJKL} applies with $\tilde\varphi=f_\C$.
By \cite[Theorem 3.9, (3.10) and \S4]{BJKL}, $\FPP(\operatorname{pgIMG}(f_p)/\Fp_p)=\FPP(\operatorname{pgIMG}(f_\C)/\C)=\FPP(\operatorname{IMG}(f_\C))$, and the latter equals $1/8$ by the computation above. For the same fixed split $p$ and the same limit inferior over $n$, the map $f_p$ is a quadratic Latt\`es map in the class covered by \cite[Theorem 2.5]{BJKL}, which gives $\liminf\ge1/8$. Hence the limit inferior equals $1/8$.
\end{proof}

\subsection{The exceptional case}

Following \cite{DH}, a rational map is exceptional if it is postcritically finite with an orbifold of finite type (spherical or Euclidean, i.e.\ of orbifold Euler characteristic $\chi\ge0$). Over $\C$ these are exactly the maps linearly conjugate to $x^{\pm d}$, to a Chebyshev polynomial $\pm T_d$, or to a Latt\`es map \cite[{\S}9, in particular Propositions 9.2 and 9.3]{DH}. For martingale groups, tame rational functions and dynamical pullbacks we use the definitions of \cite{R}.

The conjecture of Fari\~na-Asategui--Radi \cite[Conjecture~7]{FAR}, in the form stated for rational functions in \cite[Conjecture~8]{R}, concerns tamely ramified rational functions. For $f$ of degree $d\ge2$ over a field $K$ such that $G_\infty(K^{\mathrm{sep}},f,t)$ is a martingale group, it predicts that $\FPP(G_\infty)>0$ if and only if $f$ is exceptional and not linearly conjugate to $x^{\pm d}$.

\begin{corollary}\label{cor:conj}
Over a separably closed field of characteristic $0$, the prediction of \cite[Conjecture 8]{R} holds for every exceptional map, and for every tame rational function whose geometric iterated Galois group is a martingale group and which is not a dynamical pullback. For polynomials the prediction is a theorem of Fari\~na-Asategui and Radi \cite[Theorem 2]{FAR}, so the only remaining case is that of a non-exceptional dynamical pullback which is not a polynomial.
\end{corollary}

\begin{proof}
Let $f$ be tame; for Case~3 below assume in addition that $G_\infty(K^{\mathrm{sep}},f,t)$ is a martingale group (the other cases do not use this hypothesis). Replacing $K$ by $\C$ does not change the geometric group: for each $n$ the splitting field of $f^n(x)-t$ is defined over a finitely generated subfield of $K$, and embedding that subfield into $\C$ identifies the corresponding finite Galois groups (Lefschetz principle); hence the level-$n$ fixed-point counts, and therefore $\FPP$, agree. Linear conjugacy also preserves $\FPP$. Being a dynamical pullback is invariant under base change (see \cite[Definition~6.1]{R}). If $f$ is exceptional, then over $\C$ it is linearly conjugate by a M\"obius transformation to $x^{\pm d}$, to $\pm T_d$, or to a Latt\`es map \cite{DH}. 

\textbf{Case 1: $x^{\pm d}$.} The proportion is zero for $x^d$ \cite{J}. It is also zero for $x^{-d}$: for $f(z)=z^{-d}$ and $n$ odd the solutions of $f^n(z)=t$ are the roots of $z^{d^n}=t^{-1}$, while for $n$ even they are the roots of $z^{d^n}=t$; in both cases the geometric group acts by multiplication by the $d^n$-th roots of unity, only the identity fixes a root, and the level-$n$ proportion is $1/d^n\to0$. 

\textbf{Case 2: $\pm T_d$ and Latt\`es maps.} The proportion is positive for $\pm T_d$ \cite{J,FAR}, and it is positive for one-dimensional Latt\`es maps. For the last case, a Latt\`es map on $\PP^1$ has a quotient presentation $X=E/G$ with $X\cong\PP^1$; then the group $H_0=\pi_0(G)$ of linear parts is nontrivial, since a quotient of $E$ by a purely translational group is an elliptic curve and not $\PP^1$. Theorem~\ref{thm:E} therefore reduces the proportion to the fixed-origin case with nontrivial $H_0$, where the positivity criterion Corollary~\ref{cor:pos}(a) gives positivity because $h\ne1$ in $\Aut(E,0)$ has finitely many fixed points. For Latt\`es maps the martingale property follows from \cite[Theorem 5.22]{R} (see also \cite[Theorem~1.4]{HZ} for a general martingale criterion) and Theorem~\ref{thm:martingale}. 

\textbf{Case 3: not exceptional and not a dynamical pullback.} If $f$ is not exceptional and not a dynamical pullback, then \cite[Theorem 7]{R} gives $\FPP=0$. 

\textbf{Case 4: the remaining case.} The polynomial case is settled independently by \cite[Theorem 2]{FAR}, whose classification needs no martingale hypothesis, since the geometric group of every polynomial is a martingale group. Consequently the only case not covered by the results just listed is that of a non-exceptional dynamical pullback which is not a polynomial, and \cite{R} notes that there is no example of such a map whose fixed-point proportion is known.
\end{proof}

\subsection{The eight quadratic Latt\`es classes}
\label{subsec:quadratic}

Over $\C$ there are exactly eight conjugacy classes of quadratic Latt\`es maps. Their normal forms are the five expressions of \cite[Proposition~2.3]{BJKL}, with the stated multiplicities, and they are also described in \cite[{\S}8.1]{M}. Table~\ref{tab:quadratic} lists them: the unique class of orbifold type $(2,4,4)$ is the first map in \cite[Proposition~2.3]{BJKL}, the other seven have type $(2,2,2,2)$, and the associated CM fields are $\Q(i)$ (rows 1--3), $\Q(\sqrt{-2})$ (row 4) and $\Q(\sqrt{-7})$ (rows 5--8).

\begin{table}[htbp]
\centering
\small
\begin{tabular}{@{}c l l c c@{}}
\toprule
\# & Normal form (parameter) & CM field & Type & $\FPP$\\
\midrule
1 & $(x^2-2)/x^2$ & $\Q(i)$ & $(2,4,4)$ & $5/16$\\
2 & $(x^2+\alpha_1)/(x^2-\alpha_1)$, $\alpha_1=i$ & $\Q(i)$ & $(2,2,2,2)$ & $1/8$\\
3 & $(x^2+\alpha_1)/(x^2-\alpha_1)$, $\alpha_1=-i$ & $\Q(i)$ & $(2,2,2,2)$ & $1/8$\\
4 & $(x^2+\alpha_2)/(x^2-\alpha_2)$, $\alpha_2^2-2\alpha_2-1=0$ & $\Q(\sqrt{-2})$ & $(2,2,2,2)$ & $1/8$\\
5 & $(x^2+\alpha_3)/(x^2-(\alpha_3+2))$, $\alpha_3$ a root of $y^2+5y+8$ & $\Q(\sqrt{-7})$ & $(2,2,2,2)$ & $1/4$\\
6 & $(x^2+\bar\alpha_3)/(x^2-(\bar\alpha_3+2))$, $\bar\alpha_3$ the other root & $\Q(\sqrt{-7})$ & $(2,2,2,2)$ & $1/4$\\
7 & $(x^2+\alpha_3^{-1})/(x^2-(\alpha_3+2)^{-1})$, $\alpha_3$ as in row 5 & $\Q(\sqrt{-7})$ & $(2,2,2,2)$ & $1/4$\\
8 & $(x^2+\bar\alpha_3^{-1})/(x^2-(\bar\alpha_3+2)^{-1})$, $\bar\alpha_3$ as in row 6 & $\Q(\sqrt{-7})$ & $(2,2,2,2)$ & $1/4$\\
\bottomrule
\end{tabular}
\caption{The eight quadratic Latt\`es classes over $\C$ and their geometric fixed-point proportions. The normal forms are those of \cite[Proposition~2.3]{BJKL}. The first row is the unique $(2,4,4)$ class and the other seven are $(2,2,2,2)$. The values are those of Proposition~\ref{prop:quadraticFPP}.}
\label{tab:quadratic}
\end{table}

\begin{proposition}\label{prop:quadraticFPP}
Every quadratic Latt\`es map over $\C$ has geometric fixed-point proportion:
\begin{enumerate}
\item[(i)] $5/16$ for the unique class of type $(2,4,4)$;
\item[(ii)] $1/8$ for the three $(2,2,2,2)$ classes with CM field $\Q(i)$ or $\Q(\sqrt{-2})$;
\item[(iii)] $1/4$ for the four $(2,2,2,2)$ classes with CM field $\Q(\sqrt{-7})$.
\end{enumerate}
\end{proposition}
\begin{proof}
By \cite[Proposition~2.3]{BJKL} and \cite[{\S}8.1]{M}, the eight classes have quotient
presentations $E/G$ with $|G|=2$ for the seven $(2,2,2,2)$ classes and $|G|=4$ for the
$(2,4,4)$ class. In the $(2,2,2,2)$ case the group of linear parts is
$H_0=\pi_0(G)=\{\pm1\}$; for rows 2--6 $G=\{\pm1\}$, while rows 7--8 use a translated
involution. By Theorem~\ref{thm:E} (affine automorphism groups) the proportion is computed
by the determinant formula with $H_0$ in place of $G$; the only non-identity linear part is
$-1$, so it suffices to compute the dimension of the $\alpha$-active subspace $W_2$.
For the isogeny $\alpha=1\pm i$ or $\alpha=i\sqrt2$, the two eigenvalues of $\alpha$ have $2$-adic
valuation $1/2$, hence $W_2=V_2=H_1(E,\Z_2)\otimes_{\Z_2}\Q_2$ and
\[
\FPP=\frac12\abs{\det(1-(-1)\mid W_2)}_2
=\frac12\abs{4}_2=\frac18 .
\]
For the isogeny $\alpha=(1\pm i\sqrt7)/2$, the minimal polynomial is $x^2-x+2$, whose two roots have
$2$-adic valuations $0$ and $1$ (the discriminant $-7$ is a $2$-adic square); hence $W_2$
is one-dimensional, $\det(1-(-1)\mid W_2)=2$, and
\[
\FPP=\frac12\abs{2}_2=\frac14 .
\]
For the $(2,4,4)$ class, $G=\mu_4$ and the isogeny is $\alpha=1+i$; applying
Theorem~\ref{thm:main}(ii) to $\zeta=-1,\pm i$ gives
\[
\FPP=\frac14\left(\frac14+\frac12+\frac12\right)=\frac5{16}.
\]
This proves the proposition.
\end{proof}

\begin{corollary}[Arithmetic upper bounds for quadratic Latt\`es maps]\label{cor:quadraticArith}
Let $f\in F(x)$ be a quadratic Latt\`es map over a number field $F$, and let $\mathfrak p$ be a
prime of $F$ with residue field $\Fp_q$, where $q$ is odd and $q\ne7$. Assume that $f$ has good
reduction at $\mathfrak p$ and that the ramification portrait of $f$ reduces injectively modulo
$\mathfrak p$, meaning that distinct postcritical points remain distinct and the local degrees
at the postcritical points are unchanged. Then
\[
\liminf_{n\to\infty}\frac{\#\operatorname{Per}\bigl(f_{\mathfrak p},\PP^1(\Fp_{q^n})\bigr)}
{q^n+1}\le \FPP\bigl(G_\infty(\C,f_\C,t)\bigr).
\]
More precisely, the liminf is at most $1/8$ for the three $(2,2,2,2)$ classes with CM field
$\Q(i)$ or $\Q(\sqrt{-2})$, at most $1/4$ for the four $(2,2,2,2)$ classes with CM field
$\Q(\sqrt{-7})$, and at most $5/16$ for the $(2,4,4)$ class. Suppose moreover that the residue field is $\Fp_p$, where $p\equiv1\pmod4$ is a rational prime, and that $f_{\mathfrak p}$ is conjugate over $\Fp_p$ to
$(x^2+a)/(x^2-a)$ with $a^2+1=0$, i.e.\ to one of the two $\Q(i)$ normal forms of
Table~\ref{tab:quadratic} (rows 2--3). Then $\liminf=1/8$.
\end{corollary}
\begin{proof}
The restrictions $q$ odd and $q\ne7$ ensure that the reduction is tame and that the eight-class classification of Table~\ref{tab:quadratic} applies. By \cite[Corollary~3.5]{BJKL}, for the reduction $f_{\mathfrak p}$ over $\Fp_q$,
\[
\liminf_{n\to\infty}\frac{\#\operatorname{Per}(f_{\mathfrak p},\PP^1(\Fp_{q^n}))}{q^n+1}
\le \FPP\bigl(\operatorname{pgIMG}(f_{\mathfrak p})/\Fp_q\bigr).
\]
By \cite[Theorem~3.9]{BJKL} (Pink's theorem) and its proof, $f_{\mathfrak p}$ is liftable.
The injective reduction hypothesis gives $f_{\mathfrak p}$ and $f_\C$ the same ramification
portrait, and Grothendieck's specialisation theorem for tame fundamental groups (the argument
of \cite[Theorem~3.9]{BJKL}) makes their tree actions conjugate. Hence
\cite[Theorem~3.9 and (3.10)]{BJKL} applies with lift $f_\C$, and the right-hand side equals
$\FPP(\operatorname{pgIMG}(f_\C)/\C)=\FPP(G_\infty(\C,f_\C,t))$, where the last identity is
the definition of $\operatorname{pgIMG}$ over $\C$ (compare \cite[{\S}1]{BJKL}). The normal forms in rows 2--3 have the ramification portrait (2.5) of
\cite[Proposition~2.3]{BJKL}; by the injective-reduction hypothesis $f_\C$ has the same portrait,
so $f_\C$ lies in one of the classes in rows 2--4 of Table~\ref{tab:quadratic}, all with $\FPP=1/8$ by
Proposition~\ref{prop:quadraticFPP}. The final assertions follow from Proposition~\ref{prop:quadraticFPP} and \cite[Theorem~2.5]{BJKL}.
\end{proof}

\subsection{Finite-level values}

The convergence behind Theorem~\ref{thm:main} is visible at every finite level. For a scalar multiplier, Proposition~\ref{sec3:fppn-formula} gives the finite-level proportion exactly:
\[
 \FPP_n=\frac1{\abs G}\sum_{\zeta\in G}\frac{\#\bigl((1-\zeta)E[m^n]\bigr)}{\#\bigl(E[m^n]\bigr)} .
\]

The labels $(2,2,2,2)$ and $(2,4,4)$ denote the orbifold signatures of the two quotient maps of $E$ by $\mu_2$ and by $\mu_4$, respectively (cf.\ \cite{M}). For the $(2,2,2,2)$ Latt\`es map $f(x)=(x^2+1)^2/(4x(x^2-1))$ of $E:y^2=x^3-x$ with $G=\mu_2$ and $m=2$, in the numerator of the displayed formula the term for the identity is $\#\bigl((1-1)E[m^n]\bigr)=1$ and the term for $\zeta=-1$ is $\#(2E[m^n])=4^{n-1}$; hence $\FPP_n=1/8+1/(2\cdot4^n)$. The Latt\`es map changes with the curve, but the level proportions do not, so the same values of $\FPP_n$ occur for the analogous map over the non-CM curve $E:y^2=x^3+x+1$.

For the $(2,4,4)$ Latt\`es map $f(u)=(u+1)^4/(16u(u-1)^2)$ with $G=\mu_4$ and $m=2$, the term for the identity is $\#\bigl((1-1)E[m^n]\bigr)=1$, $\zeta=-1$ contributes $4^{n-1}$, and each of $\zeta=\pm i$ contributes $2^{2n-1}$; hence $\FPP_n=5/16+1/4^{n+1}$.

For $G=\mu_2$ and $m=3$ one has $\#(2E[m^n])=9^n$, hence $\FPP_n=1/2+1/(2\cdot9^n)$. Finally, for the map of \cite{BJKL}, where $G=\{\pm1\}$ and $\alpha=1+i$, the same formula with $E[\alpha^n]$ in place of $E[m^n]$ gives: the term for the identity is $\#\bigl((1-1)E[\alpha^n]\bigr)=1$ and $\zeta=-1$ contributes $\#(2E[\alpha^n])=2^{n-2}$ for $n\ge2$; thus $\FPP_n=1/8+2^{-(n+1)}$ for $n\ge2$, whereas $\FPP_1=1/2$.

\subsection{Remarks}

\begin{remark}\label{rem:covolume}
The dependence on $\alpha$ is genuine. The Haar measure is taken in the profinite Tate module $T_\alpha(E)$, so the relevant size of $(\zeta-1)T_\alpha(E)$ is its $\alpha$-adic index, not its complex covolume. The expression $1/\abs{1-\zeta}^2$ coincides with the inverse $\alpha$-adic index only when every prime divisor of $N_{F/\Q}(\zeta-1)$ divides $\deg\alpha$, where $F=\End(E)\otimes\Q$ is an imaginary quadratic field in the CM case and $F=\Q$ otherwise. For example, for $E:y^2=x^3-x$ and $G=\mu_2$ the value for $m=3$ is $1/2$, not $1/8$.
\end{remark}

\begin{remark}\label{rem:fullness-alt}
An alternative proof of the one-dimensional formula proceeds through the iterated monodromy group: \cite[Proposition~5.16 and Theorems~5.19--5.21]{R} identify the geometric iterated Galois group with the iterated monodromy group and the latter with $T_\alpha(E)\rtimes G$. The second identification is a fullness statement — the assertion that the arboreal image equals this explicit semidirect product — and it requires the cardinality hypothesis $\abs{K}\le2^{\aleph_0}$. 
\end{remark}

\section{Open problems}\label{sec:open}

\begin{enumerate}
\item \emph{Arithmetic fixed-point proportions.} All results here concern the geometric group. The arithmetic analogue, where the absolute Galois group of a number field acts, is the subject of a companion paper by the author (in preparation). Determining explicit arithmetic values in dimension at least two remains open. A concrete route is offered by the Kummer theory of Jones--Rouse \cite[Theorem~3.4]{JR}, the open-image results used there, and the bad-monodromy criterion of \cite[Definition~2.6 and Theorem~6]{FARQ8}. The companion paper proves a structure theorem for the arithmetic fixed-point proportion of Latt\`es data $(E,G,\alpha)$ over number fields. There $E$ is an elliptic curve with complex multiplication by an arbitrary order, $G\le\Aut(E)$ is a finite subgroup defined over the base field $F$, and $\alpha$ is an endomorphism in that order. For $E:y^2=x^3-x$ over $F=\Q(i)$, with $\alpha=1+i$ and $t=2+i$, it obtains the two unconditional values $31/192$ for $G=\mu_2$ and $127/384$ for $G=\mu_4$, both larger than the geometric values $1/8$ and $5/16$.
\item \emph{Quaternionic multiplication.} The determinant formula applies to abelian varieties with quaternionic multiplication (QM) through $p$-adic linear algebra (Remark~\ref{sec4:rem-qm}). To our knowledge, explicit computations for QM surfaces have not been carried out, and a classification of the finite automorphism groups and of the resulting values remains open.
\item \emph{Higher-dimensional vanishing.} We prove the vanishing criterion (Theorem~\ref{thm:zero}) in general and the T/I case distinction in the scalar case, show that there is no I-class in dimension two (Proposition~\ref{prop:dim2}), and give I-class examples in every dimension at least three. A classification of all triples $(A,G,\alpha)$ with zero fixed-point proportion is open in arbitrary dimension. In particular, to our knowledge, a higher-dimensional analogue of the classification of exceptional rational maps is not known.
\item \emph{Projective spaces.} The Latt\`es maps on $\PP^2$ are classified in \cite{Ro}. For $\PP^k$ with $k\ge3$ the classification of Latt\`es maps is incomplete. For the maps covered by the quotient definition of this paper the proportion is given by Theorem~\ref{thm:main}, but, to our knowledge, no explicit values have been tabulated for $k\ge3$.
\item \emph{Symmetrised maps.} Let $\operatorname{Sym}^k f$ be the $k$-th symmetric power of $f$. For a Latt\`es map $f$, $\operatorname{Sym}^k f$ is the Latt\`es map attached to $(A^k,G^k\rtimes S_k,\alpha^{\times k})$, so its fixed-point proportion is given by Theorem~\ref{thm:main}. Averaging over $G^k\rtimes S_k$ by cycle type (the average over each copy of $G$ contributes a factor $F$ for every cycle) gives $\FPP(\operatorname{Sym}^k f)=\frac1{k!}\sum_{\sigma\in S_k}F^{c(\sigma)}=\frac1{k!}\sum_{j=1}^k c(k,j)F^j$, where $F=\FPP(f)$, $c(\sigma)$ is the number of cycles of $\sigma$ and $c(k,j)$ is the unsigned Stirling number of the first kind (for $k=2$ this is $(F+F^2)/2$). It is natural to ask whether $\FPP(\operatorname{Sym}^k f)$ is determined by $\FPP(f)$ for an arbitrary rational map $f$.
\item \emph{Dynamical pullbacks.} The counting formula of Theorem~\ref{thm:main} is the geometric counterpart of the arithmetic formula of the companion paper, which measures the affine maps in the Galois image that have a fixed point. For a non-exceptional dynamical pullback whose defining affine map has a nonzero translation part, a fixed vertex would require an event of measure zero. We conjecture that the fixed-point proportion vanishes for such pullbacks. This is the case left open in \cite{R}.
\item \emph{Rationality.} Every value in Corollary~\ref{cor:values1} is rational. More generally Theorem~\ref{thm:main} shows that the fixed-point proportion of every Latt\`es map is rational. Is there a common derivation of the fixed-point proportion of all exceptional maps? For $\pm T_d$ the values $1/2$ ($d$ odd) and $1/4$ ($d$ even) are known \cite{J,FAR}, but a common derivation would require a better understanding of the arboreal action of Chebyshev polynomials.
\end{enumerate}

The arithmetic analogue (item 1) is the direction in which the method of this paper extends most directly.

\appendix
\section{The local index lemma}\label{app:local}

\begin{proof}[Proof of Lemma~\ref{sec4:index}]
Let $\chi(t)=\det(t-M)\in\Z_p[t]$. Write $\chi=\chi_+\chi_0$, where $\chi_+$ is the product, with multiplicities, of the irreducible factors of $\chi$ whose roots have positive $p$-adic valuation, and $\chi_0$ is the product, with multiplicities, of those whose roots are $p$-adic units. (All roots of an irreducible factor of $\chi$ over $\Q_p$ have the same $p$-adic valuation, by the argument of Lemma~\ref{sec4:val}.) Both are monic in $\Z_p[t]$. Modulo $p$ we have $\chi_+\equiv t^{\deg\chi_+}$ and $\chi_0(0)\not\equiv0$, so $\chi_+$ and $\chi_0$ are coprime in $\Fp_p[t]$. The resultant of $\chi_+$ and $\chi_0$ is congruent to
$\chi_0(0)^{\deg\chi_+}$ modulo $p$, hence is a $p$-adic unit; therefore the
Bézout relation can be chosen with $a,b\in\Z_p[t]$, that is,
\[
a\chi_+ + b\chi_0 = 1 .
\]
Put
\[
E_W:=b(M)\chi_0(M),\qquad E_U:=a(M)\chi_+(M).
\]
Since $\chi_+(M)\chi_0(M)=0$, multiplying the Bézout relation by
$\chi_+$ gives $a(M)\chi_+(M)^2=\chi_+(M)$ and hence $E_U^2=E_U$; multiplying by $\chi_0$ gives $E_W^2=E_W$ in the same way. Also $E_W+E_U=1$ and $E_WE_U=0$.
Since $\chi(M)=\chi_+(M)\chi_0(M)=0$, the two endomorphisms are complementary idempotents and
\[
\Lambda=E_W\Lambda\oplus E_U\Lambda .
\]
The summands are $M$-stable and, because $h$ commutes with $M$, also $h$-stable. They are free $\Z_p$-modules of ranks $\dim W$ and $\dim U$; hence $E_W\Lambda$ is a lattice in $W$ and $E_U\Lambda$ is a lattice in the inactive subspace $U$. Set $\Lambda_W:=E_W\Lambda$ and $\Lambda_U:=E_U\Lambda$.

On $\Lambda_U$ the endomorphism $M$ has determinant a $p$-adic unit, so $M$ is invertible and $[\Lambda_U:M\Lambda_U]=1$; hence $M^n\Lambda_U=\Lambda_U$ and
\[
[\Lambda_U:h\Lambda_U+M^n\Lambda_U]=1
\]
for every $n$.

Since $\chi_+(M)=0$ on $W$ and hence on $\Lambda_W$, and $\chi_+(M)\equiv M^{\dim W}\pmod p$, there is $N$ with $M^N\Lambda_W\subseteq p\Lambda_W$. In particular $M^n\Lambda_W\to0$ $p$-adically. We claim that
\[
\bigcap_{n\ge1}\bigl(h\Lambda_W + M^n\Lambda_W\bigr)=h\Lambda_W .
\]
The inclusion $\supseteq$ is immediate from the definition. For the reverse inclusion, let $x$ be in the intersection. For each $n$ write $x=h\lambda_n+m_n$ with $\lambda_n\in\Lambda_W$ and $m_n\in M^n\Lambda_W$. Since $m_n\to0$ and $h\Lambda_W$ is compact, a subsequence of $h\lambda_n$ converges to some $y\in h\Lambda_W$; then $x=y$, as required. This proves the claim.

Since the subgroups $h\Lambda+M^n\Lambda$ decrease with $n$, the indices $I_n:=[\Lambda:h\Lambda+M^n\Lambda]$ are non-decreasing; moreover $[\Lambda_U:h\Lambda_U+M^n\Lambda_U]=1$ for every $n$, so $I_n=[\Lambda_W:h\Lambda_W+M^n\Lambda_W]$. It therefore suffices to show that this index tends to $\abs{\det_{\Q_p}(h\mid W)}_p^{-1}$ when $\det(h\mid W)\ne0$, and to $\infty$ when $\det(h\mid W)=0$.

If $\det(h\mid W)\ne0$, then $h\Lambda_W$ is open and of finite index in $\Lambda_W$, so $\Lambda_W/h\Lambda_W$ is finite. The decreasing sequence of subgroups $(h\Lambda_W+M^n\Lambda_W)/h\Lambda_W$ of the finite group $\Lambda_W/h\Lambda_W$ has intersection $0$, hence is eventually $0$. Thus $h\Lambda_W+M^n\Lambda_W=h\Lambda_W$ for all large $n$, and
\[
\lim_{n\to\infty}[\Lambda_W:h\Lambda_W+M^n\Lambda_W]=[\Lambda_W:h\Lambda_W]
=\abs{\det_{\Q_p}(h\mid W)}_p^{-1}.
\]

If $\det(h\mid W)=0$, choose a nonzero $\Q_p$-linear functional $\varphi\colon W\to\Q_p$ vanishing on $hW$. Since $\varphi\ne0$ and $\Lambda_W$ is a lattice in $W$, the image $\varphi(\Lambda_W)$ is a nonzero fractional ideal of $\Z_p$; put $a:=v_p(\varphi(\Lambda_W))\in\Z$ and replace $\varphi$ by $p^{-a}\varphi$, so that $\varphi(\Lambda_W)=\Z_p$ (and still $\varphi(hW)=0$). By the choice of $N$ we have $M^N\Lambda_W\subseteq p\Lambda_W$, hence
\[
M^n\Lambda_W\subseteq p^{\lfloor n/N\rfloor}\Lambda_W\qquad(n\ge0),
\]
by writing $n=qN+r$ with $0\le r<N$ and using $M^r\Lambda_W\subseteq\Lambda_W$. Since $h\Lambda_W\subseteq hW$ and $\varphi(hW)=0$, it follows that
\[
\varphi\bigl(h\Lambda_W+M^n\Lambda_W\bigr)
=\varphi\bigl(M^n\Lambda_W\bigr)
\subseteq p^{\lfloor n/N\rfloor}\Z_p .
\]
Therefore $h\Lambda_W+M^n\Lambda_W\subseteq\varphi^{-1}\bigl(p^{\lfloor n/N\rfloor}\Z_p\bigr)$, and since $\varphi(\Lambda_W)=\Z_p$ the latter has index $p^{\lfloor n/N\rfloor}$ in $\Lambda_W$. Hence
\[
[\Lambda_W:h\Lambda_W+M^n\Lambda_W]\ge p^{\lfloor n/N\rfloor}\longrightarrow\infty
\]
as $n\to\infty$.

Finally, the lattices $h\Lambda$ and $M^n\Lambda$ decompose according to the direct sum, so
\[
I_n=[\Lambda:h\Lambda+M^n\Lambda]
=[\Lambda_W:h\Lambda_W+M^n\Lambda_W]\cdot
[\Lambda_U:h\Lambda_U+M^n\Lambda_U].
\]
Passing to the limit gives the statement of Lemma~\ref{sec4:index}.
\end{proof}

\section{Normalisation and the noncentral case}\label{app:twist}

This appendix records why the counting form of Theorem~\ref{thm:main}(i) needs no centrality assumption: it spells out the normalisation implicit in Proposition~\ref{sec3:fppn-formula} when $\alpha$ only normalises $G$. Condition \eqref{eq:normalize} ensures that $\alpha$ descends to the quotient $X=A/G$; in the fixed-origin case the condition is also necessary. For $\zeta\in G$ and $n\ge1$ put
\[
\gamma_n:=\alpha^{-n}\zeta\alpha^n\in G .
\]
The map $\zeta\mapsto\gamma_n$ is a bijection of the finite group $G$, with inverse given by conjugation by $\alpha^n$. In the absolute parametrisation of Lemma~\ref{sec3:affine}, an element $\sigma$ with $\sigma\mid K(A)=\zeta$ acts on $y_n$ as $\sigma(y_n)=\gamma_n y_n+c$ with $c\in A[\alpha^n]$. Thus the linear parts that occur in the counting formula are exactly the conjugates $\gamma_n$ of elements of $G$, and the average in Proposition~\ref{sec3:fppn-formula} is unchanged by the substitution $\zeta\mapsto\gamma_n$. The counting form of Theorem~\ref{thm:main}(i) therefore requires no centrality. The determinant form and the limit law require no centrality either: the determinant form follows from the proof of Theorem~\ref{thm:main}(ii) in Section~\ref{sec:det}, using Lemma~\ref{sec2:beta}, and the limit law from Remark~\ref{rem:LL-scope}.

\section*{Statements and Declarations}
\noindent\textbf{Competing interests.} The author declares no competing interests.
\noindent\textbf{Funding.} No funding was received for this work.
\noindent\textbf{Data availability.} No data sets were generated or analysed in this study; all tabulated values are derived in the text.
\noindent\textbf{Use of AI.} During the preparation of this work, the author used DeepSeek for identifying key literature, for exploring proofs, and for language editing of the manuscript. The author reviewed and verified every mathematical statement and computation independently and takes full responsibility for the content of this article. No AI tool was used as an author.

\end{document}